\documentclass{amsart}
\usepackage{amsmath,amsfonts,amssymb,amsthm,amscd, mathrsfs}
\usepackage{graphics}
\usepackage{ae, aecompl, aeguill}
\usepackage{url}
\usepackage{hyperref}

\newtheorem{thm}[equation]{Theorem}
\newtheorem{lem}[equation]{Lemma}
\newtheorem{pro}[equation]{Proposition}
\newtheorem{cor}[equation]{Corollary}
\theoremstyle{definition}
\newtheorem{defi}[equation]{Definition}
\newtheorem{nota}[equation]{Notation}
\newtheorem{rem}[equation]{Remark}
\newtheorem{exa}[equation]{Example}

\DeclareMathOperator\Supp{Supp}
\DeclareMathOperator\Aut{Aut}
\DeclareMathOperator\Bin{Bin}
\DeclareMathOperator\GL{GL}
\DeclareMathOperator\dep{dp}
\DeclareMathOperator\gr{gr}
\DeclareMathOperator\Image{Im}
\DeclareMathOperator\Id{Id}

\newcommand{\N}{\mathbb{N}}
\newcommand{\F}{\mathscr{F}}
\newcommand{\Z}{\mathbb{Z}}
\newcommand{\R}{\mathbb{R}}

\newcommand{\K}{\mathbb{K}}
\newcommand{\EL}{\mathscr{L}}
\newcommand{\gras}[1]{\textrm{\boldmath $#1$}}
\newcommand{\ie}{{\it i.e. }}

\newcommand{\cf}{{\it cf. }}
\def\G{\Gamma }
\def\D{\Delta }

\def\a{\alpha }
\def\b{\beta }
\def\g{\gamma }
\def\del{\delta }

\def\p{\varphi }
\def\m{\mathfrak }
\def\bt{\text{\scriptsize $\boxtimes$}}

\begin{document}
\author{Anatole Castella}
\title[On (twisted) LK-representations]{On (twisted) Lawrence-Krammer representations}
\date{March 20, 2008}
\maketitle

%\vspace*{-8mm}

\begin{abstract}
LK-representations (LK for Lawrence-Krammer) are linear representations of Artin-Tits monoids and groups of small type, which are of particular interest since they are known to be faithful for the monoids, and for the groups when the type is spherical, under some (strong) conditions on the defining ring $\m R$. 

For a fixed small type, they are parametrized by an $\m R$-module $\F$ for the monoid, and by a subset $\F_{\gr}$ of $\F\setminus\{0\}$ for the group. It is known that $\F_{\gr}$ is non-empty if the type has no triangle, \ie no subgraph of affine type $\tilde A_2$, and more precisely that $\F = \m R$ and $\F_{\gr} = \m R^\times$ if the type is spherical. Moreover, faithful ``twisted'' LK-representations for the non-small (spherical) crystallographic types have been constructed in [Digne, \emph{On the linearity of Artin Braid groups.} J. Algebra \textbf{268}, (2003) 39-57].

The first aim of this paper is to explicit $\F$ and $\F_{\gr}$ for any affine and small type : we establish that $\F = \m R^\N$ and $\F_{\gr} = \m R^\times \times \m R^{\N_{\geqslant 1}}$, and since this holds in particular for $\tilde A_2$, this shows that $\F_{\gr}$ can be non-empty for a graph with triangles. The second aim is to generalize the construction of \emph{op.cit.} in order to provide faithful twisted LK-representations for any Artin-Tits monoid that appears as the submonoid of fixed elements of an Artin-Tits monoid of small type under a group of graph automorphisms ; in particular, we thus get three twisted LK-representations for the Coxeter type $B_n$ (among which the one constructed in \emph{op.cit.}) and we finally show that they are pairwise non-equivalent, at least for the main choice of $\m R$. 
\end{abstract}

\section*{Introduction}

In the early 2000's, Krammer defined by explicit formulas a linear representation of the braid group on a free $\Z[x^{\pm 1},y^{\pm 1}]$-module of dimension the number of positive roots in the associated root system, and proved its faithfulness \cite{K,K2} (see also \cite{B}). This construction and the proof of faithfulness have been generalized by Cohen and Wales, and independently by Digne, to the Artin-Tits groups of spherical and small type \cite{CW,Di}, and then to all the Artin-Tits monoids of small type by Paris \cite{P} (see also \cite{H} for a short proof of the faithfulness). Since an homological version of this representation first appears in the work of Lawrence \cite{L}, those representations are commonly called \emph{Lawrence-Krammer} (LK for short). 

\medskip

In fact, LK-representations can be defined over an arbitrary unitary commutative ring $\m R$, as far as we do not require their faithfulness (see subsection \ref{proof of Hee} below for a faithfulness criterion on $\m R$). This slight generalization will give us some more insight on what is really needed in the construction, and will at least simplify some computations. Let us be more specific by considering a Coxeter graph of small type $\G$ with vertex set $I$, and its associated Artin-Tits group $B$, Artin-Tits monoid $B^+$, and set of positive roots $\Phi^+$ (see section \ref{Preliminaries}). 

Let $V$ be a free $\m R$-module with basis $(e_\a)_{\a\in \Phi^+}$. The LK-representations of $B^+$ on $V$ are linear representations $\psi : B^+\to \EL(V)$ parametrized by three elements $\m{b,\, c,\, d}$ of the unitary group $\m R^\times$ of $\m R$, and by a family $(f_i)_{i\in I}$ of linear forms on $V$ submitted to some extra conditions (see definition \ref{def lkp}). For a fixed choice of $(\m{b,c,d})\in (\m R^\times)^3$, the suitable families $(f_i)_{i\in I}$ --- that we call \emph{LK-families} --- form a submodule $\F$ of the $\m R$-module $(V^\star)^I$. When the images of $\psi$ are all invertible, then $\psi$ induces a linear representation (also called LK) $\psi_{\gr} : B\to \GL(V)$ of the Artin-Tits group $B$ ; this is precisely the case when the elements $f_i(e_{\a_i})$, $i \in I$, all belong to $\m R^\times$, and we denote by $\F_{\gr} \subseteq \F \setminus \{0\}$ the subset of LK-families that satisfy this additional condition. 

\medskip

Hence the classification of the LK-representations of $B^+$ reduces to the description of the $\m R$-module $\F$, and the question of the existence of LK-representations of $B$ reduces to the question of the non-emptiness of $\F_{\gr}$. Note moreover that, when $\G$ is connected, the elements $f_i(e_{\a_i})$, $i\in I$, are necessarily all equal for an LK-family $(f_i)_{i\in I}$. The studies of \cite{CW,Di} essentially show that, in the connected and spherical cases, an LK-family is entirely determined by the common value $\m f \in \m R$ of the $f_i(e_{\a_i})$, $i\in I$, and that this common value can be chosen arbitrarily ; hence in these cases, $\F$ is isomorphic to $\m R$, via $(f_i)_{i\in I} \mapsto f_{i_0}(e_{\a_{i_0}})$ for some $i_0 \in I$, and $\F_{\gr}$ corresponds to $\m R^\times$ via this isomorphism (see subsection \ref{cas spherique} below). This situation is partially generalized in \cite{P} where it is shown that $\F_{\gr}$ is non-empty when $\G$ has \emph{no triangle}, \ie no subgraph of affine type $\tilde A_2$ (see subsection \ref{The construction of Paris} below). As far as I know, the structure of $\F$ is not understood in general, and the question the non-emptiness of $\F_{\gr}$ is still open when $\G$ has a triangle.

\medskip

Another topic on this subject is the question of the existence of similar faithful representations in the non-small cases. A first answer is provided by \cite{Di}, where is constructed a faithful ``twisted'' LK-representation for an Artin-Tits group of type $B_n$, $F_4$ or $G_2$, using the fact that it appears as the subgroup of fixed elements under a graph automorphism, of an Artin-Tits group of type $A_{2n-1}$, $E_6$ or $D_4$ respectively. 

\medskip

The aim of this paper is to go further on those two questions. 

We first investigate the structures of $\F$ and $\F_{\gr}$ when $\G$ is of affine and small type. We show that in these cases, an LK-family $(f_i)_{i\in I}$ is not determined by the common value $\m f$ of the $f_i(e_{\a_i})$, $i\in I$, but by an infinite family $(\m f_n)_{n\in \N} \in \m R^\N$ with $\m f_0 = \m f$, which can be chosen arbitrarily ; hence $\F$ is isomorphic to $\m R^\N$, via $(f_i)_{i\in I} \mapsto (\m f_n)_{n\in \N}$, and $\F_{\gr}$ corresponds to $\m R^\times \times \m R^{\N_{\geqslant 1}}$ via this isomorphism. In particular, this holds for $\tilde A_2$ and we thus get that $\F_{\gr}$ can be non-empty when $\G$ has triangles. 

We then generalize the construction of (faithful) twisted LK-representations of \cite{Di} to any Artin-Tits monoid that appears as the submonoid of fixed elements, under a group of graph automorphisms, of an Artin-Tits monoid of small type. Note that our proof of faithfulness is different from the one of \cite{Di} as it does not use any case-by-case consideration : we show in general that the faithfulness criterion used in the small type cases also works in the twisted cases. In particular, starting with types $A_{2n}$ and $D_{n+1}$, we get two new (faithful) LK-representations of the Artin-tits group of type $B_n$. By computing the formulas obtained for a twisted LK-representation when the considered group of graph automorphisms is of order two, we show that the three twisted LK-representations for the type $B_n$ are pairwise non-equivalent, at least when $\m R = \Z[x^{\pm 1},y^{\pm 1}]$ and for the main choices of parameters.

\medskip

The paper is organized as follows.

We recall the basic results needed on Coxeter groups, root systems and Artin-Tits monoids and groups in the first section. 

In the second section, we define the LK-representations over an arbitrary commutative ring $\m R$ of Artin-Tits monoids and groups of small type following \cite{K2, CW, Di, P} (in subsection \ref{Definition}). We then generalize to our settings the faithfulness criterion on $\m R$ used in those articles and its short proof given in \cite{H} (in subsection \ref{proof of Hee}), and apply it in subsection \ref{Comments}. 

In the third section, we investigate the module of LK-families $\F$ and its subset $\F_{\gr}$ for a fixed Coxeter graph of small type and a fixed choice of parameters $\m{(b,c,d)} \in (\m R^\times)^3$. The elements of $\F$ are characterized in subsection \ref{LK-families}, and we recall the results of \cite{CW, Di} and \cite{P} on $\F$ and its subset $\F_{\gr}$ in subsections \ref{cas spherique} and \ref{The construction of Paris} respectively. Subsection \ref{The affine case} is devoted to our study of the affine case.

Finally in section \ref{Action of graph automorphisms}, we investigate the ``twisted'' LK-representations. We generalize the construction of \cite{Di} in subsection \ref{def twisted LK-reps}, and prove our ``twisted'' faithfulness criterion in subsection \ref{Faithfulness of the induced representation}. We explicit the formulas of a twisted LK-representations when the considered group of graph automorphisms is of order two in subsection \ref{formules G = 2} and apply this in subsection \ref{Twisted LK-representations of type B} to compare the twisted LK-representations of type $B_n$. As a conclusion, we give in subsection \ref{Final remark} some limitations of our approach in the non-spherical cases compared with the spherical cases of \cite{Di}.

\section{Preliminaries}\label{Preliminaries}

\subsection{General notations and definitions}\mbox{}\medskip

In all this paper, the rings we consider will be unitary, with identity element denoted by $1$ or $\Id$. Let $\m R$ be a commutative ring. We denote by $\m R^\times$ the group of units of $\m R$. If $V$ and $V'$ are two $\m R$-modules, we denote by $\mathscr{L}(V,V')$ the $\m R$-module of linear maps from $V$ to $V'$. If $V = V'$, we simply denote by $\mathscr{L}(V) = \mathscr{L}(V,V)$ the $\m R$-module of endomorphisms of $V$, by $\GL(V)$ the group of linear automorphisms of $V$ and by $V^\star = {\mathscr L}(V,\m R)$ the dual of $V$.

\medskip

A monoid is a non-empty set endowed with an associative binary operation with an identity element. A monoid $M$ is said to be \emph{left cancellative} if for any $a,\,b,\,c\in M$, $ab = ac$ implies $b = c$. The notion of \emph{right cancellativity} is defined symmetrically, and $M$ is simply said to be \emph{cancellative} when it is left and right cancellative. We denote by $\preccurlyeq$ the (left) {\it divisibility} in a monoid $M$, \ie for $a, \,b \in M$, we write $b \preccurlyeq a$ if there exists $c\in M$ such that $a = bc$ ; this leads to the natural notions of (left) gcd's and (right) lcm's in $M$. 

By \emph{a linear representation} of a monoid $M$ on an $\m R$-module $V$, we mean a monoid homomorphism $\p : M \to \mathscr{L}(V)$ ; for sake of brevity in this paper, we will often denote by $\p_b$ the image $\p(b)$ of a given $b \in M$ by a linear representation $\p$. Two linear representations $\p$ and $\p'$ of a monoid $M$, on $\m R$-modules $V$ and $V'$ respectively, are said to be equivalent if there exists a linear isomorphism $\nu : V \to V'$ such that, for every $b \in M$, $\p'_b = \nu(\p_b)\nu^{-1}$.

\subsection{Coxeter groups and Artin-Tits monoids and groups}\label{Coxeter groups and Artin-Tits monoids and groups}\mbox{}\medskip

Let $\G = (m_{i,j})_{i,j\in I}$ be a {\it Coxeter matrix}, \ie with $m_{i,j} = m_{j,i} \in \N_{\geqslant 1}\cup\{\infty\}$ and $m_{i,j} = 1 \Leftrightarrow i = j$. We will always assume in this paper that $I$ is finite ; this condition could be removed at a cost of some refinements in certain statements below (see \cite[Ch. 11]{Ca2} for some of them), which are left to the reader.

As usual, we encode the data of $\G$ by its {\it Coxeter graph}, \ie the graph with vertex set $I$, an edge between the vertices $i$ and $j$ if $m_{i,j} \geqslant 3$, and a label $m_{i,j}$ on that edge when $m_{i,j} \geqslant 4$. In the remainder of the paper, we will identify a Coxeter matrix with its Coxeter graph.

\medskip

We denote by $W = W_\G$ (resp. $B = B_\G$, resp. $B^+ = B^+_\G$) the {\it Coxeter group} (resp. {\it Artin-Tits group}, resp. {\it Artin-Tits monoid}) associated with~$\G$ : 
\begin{eqnarray*}
W & = & \langle\ s_i,\, i\in I\mid \underbrace{s_i s_j s_i \cdots}_{m_{i,j}\, \text{terms}} = \underbrace{s_j s_i s_j \cdots}_{m_{i,j}\, \text{terms}} \ \text{   if   }\  m_{i,j}\neq \infty, \text{   and   } s_i^2 = 1 \ \rangle,\\
B & = & \langle\ \gras s_i,\, i\in I\mid \underbrace{\gras s_i\gras s_j\gras s_i \cdots}_{m_{i,j}\, \text{terms}} = \underbrace{\gras s_j\gras s_i\gras s_j \cdots}_{m_{i,j}\, \text{terms}} \  \text{   if   }\ m_{i,j}\neq \infty\ \rangle,\\
B^+ & = & \langle\ \gras s_i,\, i\in I\mid \underbrace{\gras s_i\gras s_j\gras s_i \cdots}_{m_{i,j}\, \text{terms}} = \underbrace{\gras s_j\gras s_i\gras s_j \cdots}_{m_{i,j}\, \text{terms}} \  \text{   if   }\ m_{i,j}\neq \infty\ \rangle^+.
\end{eqnarray*}

Note that there is no ambiguity in writing with the same symbols the generators of $B$ and of $B^+$ since the canonical morphism $\iota : B^+ \to B$, given by the universal properties of the presentations, is injective \cite{P}, so $B^+$ can be identified with the submonoid of $B$ generated by the $\gras s_i$, $i\in I$. We denote by $\ell$ the length function on $B^+$ relatively to its generating set $\{\gras s_i \mid i\in I\}$. 

\medskip

Let $J$ be a subset of $I$. We denote by 
\begin{itemize}
\item $\G_J = (m_{i,j})_{i,j\in J}$ the submatrix of $\G$ of index set $J$,
\item $W_J = \langle s_j, \, j\in J\rangle$ the subgroup of $W$ generated by the $s_j$, $j\in J$,
\item $B_J = \langle \gras s_j, \, j\in J\rangle$ the subgroup of $B$ generated by the $\gras s_j$, $j\in J$,
\item $B^+_J = \langle \gras s_j, \, j\in J\rangle$ the submonoid of $B^+$ generated by the $\gras s_j$, $j\in J$.
\end{itemize}

It is known that $W_J$, (resp. $B_J$, resp. $B^+_J$) is the Coxeter group (resp. Artin-Tits group, resp. Artin-Tits monoid) associated with $\G_J$ (see \cite[Ch.~IV, n$^\circ$~1.8, Thm.~2]{Bo} for the Coxeter case, \cite[Ch.~II, Thm.~4.13]{vdL} for the Artin-Tits group case, the Artin-Tits monoid case being obvious). 

\medskip

We say that $J$ and $\G_J$ are {\it spherical} if $W_J$ is finite, or, equivalently, if the elements $\gras s_j$, $j \in J$, have a common (right) multiple in $B^+$. In that case, the elements $\gras s_j$, $j \in J$, have a unique (right) lcm in $B^+$, denoted by $\D_J$ and called {\it the Garside element} of $B^+_J$. Moreover, the group $B_J$ is then the {\emph group of (left) fractions} of $B^+_J$, \ie every $b \in B_J$ can be written $b = b'^{-1}b''$ with $b',\, b''\in B^+_J$ (see \cite[Props.~4.1, 5.5 and Thm.~5.6]{BS}). 

\medskip

For $b \in B^+$, we set $I(b) = \{i\in I \mid \gras s_i \preccurlyeq b\}$. In view of what has just been said, $I(b)$ is a spherical subset of $I$. 

\medskip

Let us conclude this subsection by the following easy, but fundamental, lemma : 

\begin{lem}\label{inject monoid et group}Consider a monoid homomorphism $\psi : B^+ \to G$, where $G$ is a group. Then $\psi$ extends to a group homomorphism $\psi_{\gr} : B\to G$ such that $\psi = \psi_{\gr} \circ \iota$.

Moreover if $\G$ is spherical and if $\psi$ is injective, then $\psi_{\gr}$ is injective.
\end{lem}
\proof The universal property of $B$ gives the first part. For the second, take $b \in \ker(\psi_{\gr})$ and consider a decomposition $b = b'^{-1}b''$ with $b',\, b''\in B^+$. Then $\psi_{\gr}(b) = 1$  means $\psi(b') = \psi_{\gr}(b') = \psi_{\gr}(b'') = \psi(b'')$, whence $b' = b''$ by injectivity of $\psi$ and hence $b = 1$.  \qed

\medskip

Note that if one is able to construct an injective morphism $\psi : B^+ \to G$ where $G$ is a group, then one gets that the canonical morphism $\iota$ is injective ; this is the idea of \cite{P}. In this paper, we will be interested in representations $\psi$ of $B^+$ in some linear group $\GL(V)$, hence proving their faithfulness will prove at the same time the faithfulness of the corresponding linear representation $\psi_{\gr} : B\to \GL(V)$ when $\G$ is spherical. 

\subsection{Standard root systems}\mbox{}\medskip

Let $\G = (m_{i,j})_{i,j\in I}$ be a Coxeter matrix. Details on the notions introduced here can be found in \cite{Deo}. 

\medskip

Let $E = \oplus_{i\in I}\R\a_i$ be a $\R$-vector space with basis $(\a_i)_{i\in I}$ indexed by $I$. We endow $E$ with a symmetric bilinear form $(\,.\,|\, .\,) = (\,.\,|\, .\,)_\G$ given on the basis $(\a_i)_{i\in I}$ by $(\a_i| \a_j) = -2\cos\big(\frac{\pi}{m_{i,j}}\big)$. The Coxeter group $W = W_\G$ acts on $E$ via $s_i(\b) = \b-(\b|\a_i)\a_i$. 

The ({\it standard}) {\it root system} associated with $\G$ is by definition the set $\Phi = \Phi_\G = \{w(\a_i) \mid w\in W, \, i\in I\}$. It is well-known that $\Phi = \Phi^+ \sqcup \Phi^-$, where $\Phi^+ = \Phi \bigcap \left(\oplus_{i\in I}\R^+\a_i\right)$ and $\Phi^- = -\Phi^+$. For $\a = \sum_{i\in I}\lambda_i\a_i \in \Phi$ we call \emph{support} of $\a$ the set $\Supp(\a) = \{i\in I \mid \lambda_i \neq 0\}$.

\medskip

We will always represent a subset $\Psi$ of $\Phi^+$ by a graph with vertex set $\Psi$ and an edge labeled $i$ between two vertices $\a$ and $\b$ if $\a = s_i(\b)$. For example, the situation where $\b$ is fixed by $s_i$ will be drawn by a loop
 \begin{picture}(35,7)(5,0)
\put(20,3){\circle*{4}}
\put(25,3){\circle{10}}
\put(32,0){\footnotesize $i$}
\put(8,0){$\b$}
\end{picture}.

\medskip

Such a graph is naturally $\N$-graded via the {\it depth} function on $\Phi^+$, where the depth of a root $\a \in \Phi^+$ is by definition $\dep(\a) = \min\{l(w) \mid w \in W, \, w(\a) \in \Phi^-\}$. Contrary to what suggests this terminology, in all the graphs that we will draw, we chose to place a root of great depth {\it above} a root of small depth ; so drawings like the following ones (with $\b$ {\it above} $\a$), will all mean that $\b = s_i(\a)$ (or equivalently $\a = s_i(\b)$) and $\dep(\b) > \dep(\a)$ :

\begin{center}
\begin{picture}(30,25)
 \put(0,0){\circle*{4}}\put(5,-3){$\a$}
 \put(13,20){\circle*{4}}\put(18,20){$\b$}
 \put(0,0){\line(2,3){12}}\put(0,8){\footnotesize $i$}
\end{picture},
\begin{picture}(30,25)
 \put(10,0){\circle*{4}}\put(15,-3){$\a$}
 \put(10,20){\circle*{4}}\put(15,20){$\b$}
 \put(10,0){\line(0,1){20}}\put(5,8){\footnotesize $i$}
\end{picture},
\begin{picture}(30,25)
 \put(13,0){\circle*{4}}\put(20,-3){$\a$}
 \put(0,20){\circle*{4}}\put(5,20){$\b$}
 \put(13,0){\line(-2,3){12}}\put(9,8){\footnotesize $i$}
\end{picture}
$\ldots$
\end{center}

\medskip

\begin{lem}[{\cite[Lem 1.7]{BH}}]\label{lemme BH}
 Let $i \in I$ and $\a \in \Phi^+ \setminus \{\a_i\}$. Then 
$$\dep(s_i(\a)) = 
\begin{cases}
\dep(\a)-1 & \text{if  } (\a|\a_i)>0,\\
\dep(\a) & \text{if  } (\a|\a_i)=0,\\
\dep(\a)+1 & \text{if  } (\a|\a_i)<0,\\
\end{cases}$$
\end{lem}
%\proof This is . \qed

\medskip

In the remainder of the paper, we will often consider subsets of $\Phi^+$ of the form $\{w(\a)\mid w \in W_{\{i,j\}}\} \bigcap \Phi^+$, for $\a \in \Phi^+$ and $i,\,j\in I$ with $m_{i,j} = 2$ or $3$, so the following definition and remark will be useful :

\begin{defi}Let $\a \in \Phi^+$ and $J \subseteq I$. We call \emph{$J$-mesh} of $\a$, or simply \emph{mesh}, the set $[\a]_J := \{w(\a)\mid w \in W_J\} \bigcap \Phi^+$. This terminology is inspired by personal communications with H\'ee.
\end{defi}

\begin{rem}\label{mailles possibles}Let $\a \in \Phi^+$ and $i,\,j\in I$ with $m_{i,j} = 2$ or $3$. Then, up to exchanging $i$ and $j$, the graph of the mesh $[\a]_{\{i,j\}}$ is one of the following :
\begin{itemize}
 \item if $m_{i,j} = 2$ : 
\end{itemize}
\vspace{-4mm}
\begin{center}
\begin{tabular}[t]{|c|c|c|c|} \hline  Type 1 & Type 2 & Type 3 & Type 4 \\ \hline \hline 
 \begin{picture}(40,10)(0,-7)
\put(20,3){\circle*{4}}
\put(25,3){\circle{10}}
\put(32,0){\small $j$}
\put(6,0){$\a_i$}
\end{picture} &  
 \begin{picture}(40,10)(0,-7)
\put(20,3){\circle*{4}}
\put(25,3){\circle{10}}
\put(32,0){\small $j$}
\put(15,3){\circle{10}}
\put(5,0){\small $i$}
%\put(18,10){$\a$}
\end{picture} &
\begin{picture}(40,22)(-5,0)
%\put(20,0){\circle*{4}}
\put(20,0){\circle*{4}}
\put(20,20){\circle*{4}}
\put(20,0){\line(0,1){20}}
\put(22,8){\small $i$}
\put(15,0){\circle{10}}
\put(4,-2){\small $j$}
\put(15,20){\circle{10}}
\put(4,18){\small $j$}
%\put(12,2){$\a$}
%\put(12,22){$\b$}
\end{picture} & 
\begin{picture}(50,35)(-5,10)
\put(20,0){\circle*{4}}
\put(40,20){\circle*{4}}
\put(20,40){\circle*{4}}
\put(0,20){\circle*{4}}
\put(20,0){\line(1,1){20}}
\put(20,0){\line(-1,1){20}}
\put(20,40){\line(1,-1){20}}
\put(20,40){\line(-1,-1){20}}
\put(33,5){\small $i$}
\put(33,30){\small $j$}
\put(2,5){\small $j$}
\put(5,30){\small $i$}
%\put(12,2){$\a$}
%\put(12,22){$\b$}
\end{picture} \\ & & & \\ \hline 
\end{tabular}
\end{center}
%\begin{itemize}

%\vspace{2mm}
\begin{itemize}
\item if $m_{i,j} = 3$ :
\end{itemize}
\vspace{-4mm}
\begin{center}
\begin{tabular}[t]{|c|c|c|c|}\hline Type 5 & Type 6 & Type 7 & Type 8 \\ \hline \hline 
\begin{picture}(50,25)(10,5)
%\put(20,0){\circle*{4}}
\put(55,00){\circle*{4}}
\put(35,20){\circle*{4}}
\put(15,00){\circle*{4}}
%\put(20,0){\line(1,1){20}}
%\put(20,0){\line(-1,1){20}}
\put(35,20){\line(1,-1){20}}
\put(35,20){\line(-1,-1){20}}
\put(18,10){\small $j$}
\put(48,10){\small $i$}
\put(10,-10){$\a_i$}
\put(52,-10){$\a_j$}
\put(21,25){$\a_i\!+\!\a_j$}
%\put(12,22){$\b$}
\end{picture} &
\begin{picture}(40,10)(0,-2)
\put(20,3){\circle*{4}}
\put(25,3){\circle{10}}
\put(32,0){\small $j$}
\put(15,3){\circle{10}}
\put(5,0){\small $i$}
%\put(18,10){$\a$}
\end{picture} &
\begin{picture}(40,42)(0,15)
%\put(20,0){\circle*{4}}
\put(20,0){\circle*{4}}
\put(20,40){\circle*{4}}
\put(20,20){\circle*{4}}
\put(20,0){\line(0,1){40}}
\put(22,8){\small $j$}
\put(22,28){\small $i$}
\put(15,0){\circle{10}}
\put(4,-2){\small $i$}
\put(15,40){\circle{10}}
\put(4,38){\small $j$}
%\put(12,2){$\a$}
%\put(12,22){$\b$}
\end{picture} &
\begin{picture}(50,42)(-5,23)
\put(20,0){\circle*{4}}
\put(40,20){\circle*{4}}
\put(40,40){\circle*{4}}
\put(20,60){\circle*{4}}
\put(0,20){\circle*{4}}
\put(0,40){\circle*{4}}
\put(20,0){\line(1,1){20}}
\put(0,20){\line(0,1){20}}
\put(40,20){\line(0,1){20}}
\put(20,0){\line(-1,1){20}}
\put(20,60){\line(1,-1){20}}
\put(20,60){\line(-1,-1){20}}
\put(33,5){\small $j$}
\put(33,50){\small $j$}
\put(2,5){\small $i$}
\put(5,50){\small $i$}
\put(-6,26){\small $j$}
\put(42,26){\small $i$}
%\put(12,2){$\a$}
%\put(12,22){$\b$}
\end{picture} \\ & & & \\ & & & \\ \hline 
\end{tabular}
\end{center}
%\end{itemize}
\end{rem}

\medskip

Let $J$ be a subset of $I$. We denote by $\Phi_J$ the subset $\{w(\a_j)\mid w \in W_J, \, j\in J\}$ of $\Phi$. It is clear that $\Phi_J$ is the root system associated with $\G_J$ in $\oplus_{j\in J}\R\a_j$.

\subsection{Graph automorphisms}\label{Graph automorphisms}\mbox{}\medskip

We call {\it automorphism} of a Coxeter matrix $\G = (m_{i,j})_{i,j\in I}$ every permutation $g$ of $I$ such that $m_{g(i),g(j)} = m_{i,j}$ for all $i,\,j \in I$, and we denote by $\Aut(\G)$ the group the constitute.

%\medskip

Any automorphism of $\G$ clearly acts by automorphisms on $W$, $B$ and $B^+$ by permuting the corresponding generating set. If $G$ is a subgroup of $\Aut(\G)$, we denote by $W^G$, $B^G$ and $(B^+)^G$ the corresponding subset of fixed points under the action of the elements of $G$. It is known that $W^G$ (resp. $(B^+)^G$) is a Coxeter group (resp. Artin-Tits monoid) associated with a certain Coxeter graph $\G'$ easily deduced from $\G$, and the analogue holds for $B^G$ when $\G$ is spherical, or more generally of {\it FC-type} (see \cite{H1,Mu} for the Coxeter case, \cite{Mi,C2,C2',Ca} for the Artin-Tits case). Note that the standard generator of $(B^+)^G$ are the Garside elements $\D_J$ of $B^+_J$, for $J$ running through the spherical orbits of $I$ under $G$. 

%\medskip

Similarly, any automorphism $g$ of $\G$ acts by a linear automorphism on $E = \oplus_{i\in I}\R\a_i$ by permuting the basis $(\a_i)_{i\in I}$. This action stabilizes $\Phi$ and $\Phi^+$, and the induced action on those sets is given by $w(\a_i) \mapsto (g(w))(\a_{g(i)})$. 

\section{LK-representations}\label{LK-representations}

In subsection \ref{Definition} below, we define the Lawrence-Krammer representations, over an arbitrary (unitary) commutative ring $\m R$, of the Artin-Tits monoids an groups of small type. The definition is inspired by the ones of \cite{K2, CW, Di, P}, where $\m R$ is chosen to be $\Z[x^{\pm 1},y^{\pm 1}]$ (\cf section \ref{Comments} below). 

In subsection \ref{proof of Hee}, we extend to our settings the faithfulness criterion of \cite{K2, CW, Di, P} and prove it following \cite{H}. We apply that criterion in subsection \ref{Comments}. %: under some (strong) conditions on the ring $\m R$ and on the parameters of a given Lawrence-Krammer representation of an Artin-Tits monoid (see theorem \ref{criterion} below), then it is faithful.

\medskip

From now on, we assume that $\G = (m_{i,j})_{i,j\in I}$ is a Coxeter matrix of \emph{small} type, \ie with $m_{i,j} \in \{1,2,3\}$ for all $i,\, j \in I$.

\subsection{Definition}\label{Definition}\mbox{}\medskip

Let $\mathfrak R$ be a commutative ring and $V$ be a free $\m R$-module with basis $(e_\a)_{\a\in \Phi^+}$ indexed by $\Phi^+$. 

\begin{nota}
 For $f \in V^\star$ and $e \in V$, we denote by $f\bt e$ the element of $\mathscr{L}(V)$ given by $(f\bt e)(v) = f(v)e$ for every $v \in V$.
\end{nota}

\begin{rem}\label{formulaire de calcul}
 Consider $\p \in {\mathscr L}(V)$, $f,\, f' \in V^\star$ and $e,\,e' \in V$. Then :
\begin{enumerate}
 \item $\p(f\bt e) = f\bt \p(e)$,
\item $(f\bt e) \p = (f \p)\bt e$,
\item $(f\bt e) (f'\bt e') = f(e')(f'\bt e)$.
\end{enumerate}
\end{rem}

\begin{defi}\label{def psii}Fix $\m {(a, \, b,\, c,\, d) \in R}^4$ and a family of linear forms $(f_i)_{i\in I} \in (V^\star)^I$. For $i \in I$, we denote by $\m f_{i,\a}$ the element $f_i(e_\a)$, for $\a \in \Phi^+$, and by 
\begin{itemize}
 \item $\p_i$ the endomorphism of $V$ given on the basis $(e_\a)_{\a\in \Phi^+}$ by 
\end{itemize}
\begin{center}
$\begin{cases}
 \p_i(e_{\a}) = 0 & \text{  if  }\ \ \a = \a_i, \\
\p_i(e_\a) = \m d e_\a & \text{  if  
\begin{picture}(36,12)(0,-2)
\put(20,0){\circle*{4}}
\put(25,0){\circle{10}}
\put(32,-3){\footnotesize $i$}
\put(10,2){$\a$}
\end{picture}},\\
\begin{cases}\p_i(e_\b) = \m b e_\a \\ \p_i(e_\a) = \m a e_\a + \m c e_\b \end{cases} & \text{  if  
\begin{picture}(32,20)(0,8)
%\put(20,0){\circle*{4}}
\put(20,0){\circle*{4}}
\put(20,20){\circle*{4}}
\put(20,0){\line(0,1){20}}
\put(22,8){\footnotesize $i$}
\put(10,0){$\a$}
\put(10,18){$\b$}
\end{picture} in   } \, \Phi^+.
\end{cases}$
\end{center}

\begin{itemize}
\item $\psi_i$ the endomorphism of $V$ given by $\psi_i = \p_i+f_i\bt e_{\a_i}$.
\end{itemize}
\end{defi}

\begin{rem}\label{les matrices}
 If one fixes an arrangement of the basis $(e_\a)_{\a \in \Phi^+}$ so that $e_{\a_i}$ is the leftmost element and $e_\b$ is the right successor of $e_\a$ whenever $\b = s_i(\a)$ with $\dep(\b)>\dep(\a)$, then the matrix of $\p_i$ in this basis is block diagonal, with blocks 

\begin{center}
$\left\{\begin{array}{cl}
  \overset{\text{\ \scriptsize $e_{\a}$}}{\big(\begin{matrix}0\end{matrix}\big)}
 & \text{ if }\ \a = \a_i, \\
\overset{\text{\ \scriptsize $e_{\a}$}}{\big(\begin{matrix}\m d\end{matrix}\big)}
 & \text{ if \begin{picture}(32,20)(5,-2)
\put(20,0){\circle*{4}}
\put(25,0){\circle{10}}
\put(32,-3){\footnotesize $i$}
\put(10,2){$\a$}
\end{picture}}, \\
  \overset{\text{\scriptsize \ \ $e_\a$\ \ $e_\b$}}{\begin{pmatrix}\m a & \m b \\ \m c & 0 \end{pmatrix}}
  & \text{ if \begin{picture}(25,30)(5,8)
%\put(20,0){\circle*{4}}
\put(20,0){\circle*{4}}
\put(20,20){\circle*{4}}
\put(20,0){\line(0,1){20}}
\put(22,8){\footnotesize $i$}
\put(10,0){$\a$}
\put(10,20){$\b$}
\end{picture}}.
\end{array}\right.$
\end{center}

And the matrix of $\psi_i$ is the same except that the first row (the one of index ${\a_i}$), which is zero in $\p_i$, if replaced by the row $(\m f_{i,\a})_{\a \in \Phi^+} = (f_i(e_\a))_{\a \in \Phi^+}$.
\end{rem}

\medskip

Let us now exhibit conditions on $\m{a,\, b,\, c,\, d}$ and $(f_i)_{i\in I}$ so that the map $\gras s_i \mapsto \psi_i$ extends to a linear representation $\psi$ of $B^+$ into $\EL(V)$ or $\GL(V)$.

\begin{lem}\label{condition d'inversibilité}
The map $\psi_i$ is invertible if and only if $\m b, \, \m c,\, \m d$ and $\m f_{i,\a_i}$ belong to $\m R^\times$, in which case the inverse of $\psi_i$ is given by 
\begin{center}
$\begin{cases}
 \psi^{-1}_i(e_{\a}) = \dfrac{1}{\m f_{i,\a_i}}e_{\a_i} & \text{  if  }\ \ \a = \a_i, \\
\psi^{-1}_i(e_\a) = \dfrac{1}{\m d}\Big( e_\a - \dfrac{\m f_{i,\a}}{\m{f}_{i,\a_i}}e_{\a_i}\Big) & \text{  if  
\begin{picture}(36,12)(0,-2)
\put(20,0){\circle*{4}}
\put(25,0){\circle{10}}
\put(32,-3){\footnotesize $i$}
\put(10,2){$\a$}
\end{picture}},\\
\begin{cases}\psi^{-1}_i(e_\b) = \dfrac{1}{\m c}\Big( e_\a - \dfrac{\m a}{\m{b}}e_\b + \dfrac{\m{af}_{i,\b}-\m {bf}_{i,\a}}{\m{bf}_{i,\a_i}}e_{\a_i}\Big)\\ \psi^{-1}_i(e_\a) = \dfrac{1}{\m b}\Big( e_\b - \dfrac{\m f_{i,\b}}{\m{f}_{i,\a_i}}e_{\a_i}\Big) \end{cases} & \text{  if  
\begin{picture}(32,20)(0,8)
%\put(20,0){\circle*{4}}
\put(20,0){\circle*{4}}
\put(20,20){\circle*{4}}
\put(20,0){\line(0,1){20}}
\put(22,8){\footnotesize $i$}
\put(10,0){$\a$}
\put(10,18){$\b$}
\end{picture} in   } \, \Phi^+.
\end{cases}$
\end{center}
\end{lem}
\proof Straightforward computations. \qed

\begin{lem}\label{condition sur les phii}Consider $i,\, j\in I$ with $i\neq j$. 
\begin{enumerate}
 \item If $m_{i,j} = 2$, then $\p_i\p_j = \p_j \p_i$.
\item If $m_{i,j} = 3$ and if $\mathfrak {a\big(d(a-d)+bc\big)} = 0$, then $\p_i\p_j\p_i = \p_j \p_i \p_j$.
\end{enumerate}
\end{lem}
\proof For every $\a \in \Phi^+$, the linear maps $\p_i$ and $\p_j$ stabilize the submodule of $V$ generated by the elements $e_\b$ for $\b$ running through the $\{i,j\}$-mesh $[\a]_{i,j}$ of $\a$. The results then follow from the direct computations of the matrices of the restrictions of $\p_i\p_j$ and $\p_i\p_j\p_i$ to those submodules (of dimension 1, 2, 3, 4 or 6 in view of remark \ref{mailles possibles}). Note that the only case where the condition $\mathfrak {a\big(d(a-d)+bc\big)} = 0$ is needed is the case of a mesh of type 7 in the nomenclature of remark \ref{mailles possibles}. \qed

\begin{lem}\label{condition sur les psii}Consider $i,\, j\in I$ with $i\neq j$, and assume that $\mathfrak {d(a-d)+bc} = 0$ and $f_{i}(\a_j) = f_{j}(\a_i) = 0$. 
\begin{enumerate}
 \item If $m_{i,j} = 2$, then $\psi_i\psi_j = \psi_j \psi_i$ if and only if $f_i\p_j = \m df_i$ and $f_j\p_i = \m df_j$.
\item If $m_{i,j} = 3$, then $f_i\p_j = f_j\p_i$ implies $\psi_i\psi_j\psi_i = \psi_j \psi_i \psi_j$, and the converse is true if $\m c \in \m R^\times$.
\end{enumerate}
\end{lem}
\proof Note that, since $f_{i}(\a_j) = 0$, we get, by using the formulas of remark \ref{formulaire de calcul} : $(f_i\bt v)(f\bt e_{\a_j}) = 0$ for every $(v,f) \in V\times V^\star$ (and similarly if we exchange $i$ and $j$), and hence 
\begin{itemize}
 \item $\psi_i\psi_j = \p_i\p_j + f_j\bt \p_i(e_{\a_j}) + (f_i\p_j)\bt e_{\a_i}$, and 
\item $\psi_i\psi_j\psi_i = \p_i\p_j\p_i + f_i\bt \p_i\p_j(e_{\a_i}) + (f_j\p_i)\bt \p_i(e_{\a_j})$

\hspace*{5cm} $ + \big(f_i\p_j\p_i + (f_i\p_j)(e_{\a_i})f_i\big)\bt e_{\a_i}$.
\end{itemize}
If $m_{i,j} = 2$, then $\p_i(e_{\a_j}) = \m d e_{\a_j}$, thus we get, by symmetry in $i$ and $j$ and by case (i) of the previous lemma : 
$$\psi_i\psi_j - \psi_j\psi_i = (f_i\p_j - \m d f_i)\bt e_{\a_i} - (f_j\p_i - \m d f_j)\bt e_{\a_j}.$$
This establishes (i). If $m_{i,j} = 3$, then $\p_i(e_{\a_j}) = \m a e_{\a_j} + \m c e_{\a_i+\a_j}$, $\p_i\p_j(e_{\a_i}) = \m {bc} e_{\a_j}$, thus we get, by symmetry in $i$ and $j$ and by case (ii) of the previous lemma : 
\begin{eqnarray*}
 \psi_i\psi_j\psi_i - \psi_j\psi_i\psi_j & = & \m c (f_j\p_i-f_i\p_j)\bt e_{\a_i+\a_j} \\ & & + \big(f_j\p_i(\p_j-\m a \Id) + (f_j\p_i)(e_{\a_j})f_j-\m{bc}f_i\big)\bt e_{\a_j} \\ 
& & - \underbrace{\big(f_i\p_j(\p_i-\m a \Id) + (f_i\p_j)(e_{\a_i})f_i-\m{bc}f_j\big)}_{F_{i,j}}\bt e_{\a_i}.
\end{eqnarray*}
The second part of (ii) is now clear, and to show the direct implication, we have to show that $f_i\p_j = f_j\p_i$ implies $F_{i,j} = 0$ (this will give $F_{j,i} = 0$ by symmetry and hence $\psi_i\psi_j\psi_i = \psi_j\psi_i\psi_j$). But since $\p_i(e_{\a_i}) = 0$, we have 
$$ F_{i,j} = (f_i\p_j-f_j\p_i)(\p_i-\m a \Id) + (f_i\p_j-f_j\p_i)(e_{\a_i})f_i+f_j\big(\p_i^2 - \m a \p_i -\m{bc}\Id\big),$$
and the linear form $f_j\big(\p_i^2 - \m a \p_i -\m{bc}\Id\big)$ is the zero form, since it is zero on 
\begin{itemize}
 \item $e_{\a_i}$ since $\p_i(e_{\a_i}) = 0$ and $f_j(e_{\a_i}) = 0$,
\item $e_\a$ if 
\begin{picture}(36,12)(5,-2)
\put(20,0){\circle*{4}}
\put(25,0){\circle{10}}
\put(32,-3){\footnotesize $i$}
\put(10,2){$\a$}
\end{picture}
 since then $\p_i(e_\a) = \m de_\a$ and $\m d^2 - \m{ad-bc} = 0$,
\item $e_\a$ and $e_\b$ if 
\begin{picture}(32,20)(0,8)
%\put(20,0){\circle*{4}}
\put(20,0){\circle*{4}}
\put(20,20){\circle*{4}}
\put(20,0){\line(0,1){20}}
\put(22,8){\footnotesize $i$}
\put(10,0){$\a$}
\put(10,18){$\b$}
\end{picture} in $\Phi^+$, since $X^2 - \m a X - \m{bc}$ is the characteristic 

\vspace{3mm}

\noindent polynomial of the restriction of $\p_i$ on $\m R e_\a\oplus \m R e_\b$ (see remark \ref{les matrices}).
\end{itemize}

Whence the result.\qed

\begin{defi}\label{def lkp}Fix $\m{(b,c,d)} \in (\m R^\times)^3$, set $\m a = \m d - \dfrac{\m{bc}}{\m d}$, and consider the linear maps $\p_i \in \EL(V)$, $i\in I$, as in definition \ref{def psii}. We say that a family $(f_i)_{i\in I} \in (V^\star)^I$ is an \emph{LK-family} (\emph{relatively to} $(\m{b,c,d})$) if it satisfies the following properties :
\begin{enumerate}
%\item $\m b, \, \m c,\, \m d \in \m R^\times$ and $\mathfrak {d(a-d)+bc} = 0$,
\item for $i,\,j \in I$ with $i\neq j$, $f_i(e_{\a_j})=0$,
\item for $i,\, j \in I$ with $m_{i,j} = 2$, $f_i\p_j = \m df_i$,
\item for $i,\, j \in I$ with $m_{i,j} = 3$, $f_i\p_j = f_j\p_i$.
\end{enumerate}

We denote by $\F = \F_{(\m{b,c,d})}$ the set of LK-families relatively to $(\m{b,c,d})$. This is clearly a submodule of the $\m R$-module $(V^\star)^I$. We denote by $\F_{\gr}$ the subset of $\F$ composed of the LK-families for which $f_i(e_{\a_i}) \in \m R^\times$ for every $i \in I$.

In view of lemma \ref{condition sur les psii} above, for every LK-family $(f_i)_{i\in I}$, the map $\gras s_i \mapsto \psi_i = \p_i + f_i\bt e_{\a_i}$ extends to a linear representation $\psi = \psi_{\m{(b,c,d)},(f_i)_{i\in I}} : B^+ \to \EL(V)$. Moreover, if $(f_i)_{i\in I} \in \F_{\gr}$, then in view of lemma \ref{condition d'inversibilité}, the images of $\psi$ are invertible, and hence $\psi : B^+ \to \GL(V)$ induces a linear representation $\psi_{\gr} : B \to \GL(V)$.

We call \emph{Lawrence-Krammer representation} --- \emph{LK-representation} for short --- the representation $\psi$ of $B^+$ and, when appropriate, the representation $\psi_{\gr}$ of $B$.
\end{defi}

\begin{rem}
The assumption on $\m{b,\, c,\, d}$ to be units of $\m R$ is not needed to define the LK-representations of $B^+$ and for the faithfulness criterion of the following subsection. We included it in the definition since we are mainly interested in LK-representations of $B^+$ that extends to LK-representations of $B$, and since it will be of importance in our general study of LK-families in section \ref{On LK-families} below.
\end{rem}

\subsection{Faithfulness criterion}\label{proof of Hee}\mbox{}\medskip

The key argument in \cite{K, K2, CW, Di, P} is that the LK-representation $\psi$ they consider is faithful. The faithfulness criterion used each time can be summarized as follows (where $\m{a,\, b,\, c,\, d}$, $(f_{i})_{i\in I}$ and $\psi$ are as in definition \ref{def lkp}) :

\begin{thm}\label{criterion}Assume that the following two conditions are satisfied : 
\begin{enumerate}
\item $\Image(\psi)$ is a left cancellative submonoid of $\EL(V)$, 
\item there exists a totally ordered commutative ring $\m R_0$ and a ring homomorphism $\m R \to \m R_0$, $\m x \mapsto \overline{\m x}$, such that $\overline{\m a}$, $\overline{\m b}$, $\overline{\m c}$, $\overline{\m d}$ are positive and $\overline{\m f_{i,\a}} = 0$ for every $(i,\a) \in I\times \Phi^+$. 
\end{enumerate}

Then the LK-representation $\psi$ is faithful.
\end{thm}

In the remainder of this subsection, we sketch the (much easier) proof of this criterion obtained by H\'ee in \cite{H}. It does not involve any consideration on closed sets of positive roots, nor on the maximal simple (left) divisor of an element of $B^+$, and rely only on the two following (elementary) lemmas and a look at the defining formulas of $\psi$.

\begin{lem}[{\cite[Prop. 1]{H}}]\label{Hee1}Let $\rho : B^+ \to M$ be a monoid homomorphism where $M$ is left cancellative. If $\rho$ satisfies $\rho(b) = \rho(b') \Rightarrow I(b)= I(b')$ for all $b,\,b' \in B^+$, then it is injective.
\end{lem}
\proof Under those assumptions, one can show, by induction on $\ell(b)$, that $\rho(b) = \rho(b')$ implies $b = b'$. (See lemma \ref{reduction2} below for a twisted version of that result.) \qed

\begin{nota}\label{notations bin}
 We denote by $\Bin(\Omega)$ the monoid of binary relations on a set $\Omega$, where the product $RR'$ of two binary relations $R$ and $R'$ is defined on $\Omega$ by $\b RR'\a \Leftrightarrow \exists \, \gamma \in \Omega$ such that $\b R \gamma$ and $\gamma R'\a$. 

For $R\in \Bin(\Omega)$ and $\Psi \subseteq \Omega$, we denote by $R(\Psi)$ the set $\{ \b \in \Omega \mid \exists \, \a \in \Psi, \ \b R\a\}$. 
\end{nota}

We will use again the following lemma in the proof of our faithfulness criterion in subsection \ref{Faithfulness of the induced representation} below. So, for completeness, we detail its proof here.

\begin{lem}[{\cite[Prop. 2]{H}}]\label{Hee3}Let $B^+ \to \Bin(\Omega)$, $b \mapsto R_b$, be a monoid homomorphism, and let $(\a_i)_{i\in I}$ be a family of elements of $\Omega$ such that 
\begin{enumerate}
\item $\a_i \not \in R_{\gras s_i}(\Omega)$,
\item if $i \neq j$, then $\a_i R_{\gras s_j} \a_i$,
\item if $m_{i,j} = 3$, then $\a_i R_{\gras s_j}R_{\gras s_i} \a_j$.
\end{enumerate}

Then for every $b \in B^+$, we have $\gras s_i \preccurlyeq b \Leftrightarrow \a_i \not \in R_b(\Omega)$. In particular, for $b,\, b' \in B^+$, we get $R_b(\Omega) = R_{b'}(\Omega) \Rightarrow I(b) = I(b')$.
\end{lem}
\proof Since the map $b \to R_b$ is a monoid homomorphism, we get that $b' \preccurlyeq b$ implies $R_b(\Omega) \subseteq R_{b'}(\Omega)$ for every $b$, $b' \in B^+$. Thanks to property (i), this shows that $\gras s_i \preccurlyeq b \Rightarrow \a_i \not \in R_b(\Omega)$. For the converse, assume that $\gras s_i \not \preccurlyeq b$ and let us prove by induction on $\ell(b)$ that $\a_i \in R_b(\Omega)$. If $b = 1$, then $R_b$ is the equality relation and hence $\a_i \in R_b(\Omega)$. If $\ell(b) > 0$, fix $j \in I$ such that $\gras s_j \preccurlyeq b$ (hence $j \neq i$) and denote by $b_1$ the element of $B^+$ such that $b = \gras s_jb_1$ ; if $\gras s_i \not \preccurlyeq b_1$, then $\a_i \in R_{b_1}(\Omega)$ by induction and hence $\a_i \in R_{b}(\Omega)$ thanks to property (ii) ; if $\gras s_i \preccurlyeq b_1$, then $b = \gras s_j \gras s_i b_2$ for some $b_2 \in B^+$, and since $\gras s_i \not \preccurlyeq b$, we necessarily have $m_{i,j} \neq 2$ (hence $m_{i,j} = 3$) and $\gras s_j \not \preccurlyeq b_2$, so $\a_j \in R_{b_2}(\Omega)$ by induction and $\a_i \in R_b(\Omega)$ thanks to property~(iii). \qed

\begin{rem}\label{des matrices positives aux relations}
Let $\m R_0$ be a totally ordered commutative ring and let $V_0$ be a free $\m R_0$-module with basis $(e_\a)_{\a \in \Omega}$. We denote by $\m R_0^+$ the set (semiring) of non-negative elements of $\m R_0$ and by $\mathscr{L}^+(V_0)$ the submonoid of $\mathscr{L}(V_0)$ composed of the linear maps $\p : V_0 \to V_0$ such that $\p(e_\a) \in \oplus_{\b \in \Omega}\m R_0^+e_\b$ for all $\a \in \Omega$.

Then there is a monoid homomorphism $\mathscr{L}^+(V_0) \to \Bin(\Omega)$, $\p \mapsto R_\p$, where $R_\p$ is given by 
%\begin{center}
 $\b R_\p \a \Leftrightarrow$ the coefficient of $e_\b$ in $\p(e_\a)$ is positive.
%\end{center}
%is a monoid homomorphism.  
\end{rem}

Now assume that we are in the situation of condition (ii) of theorem \ref{criterion}. 

\begin{defi}[{\cite[4.3]{H}}]\label{relations de Digne}If we denote by $V_0$ the free $\m R_0$-module with basis $(e_\a)_{\a \in \Phi^+}$, then the ring homomorphism $\m R \to \m R_0$, $\m x \mapsto \overline{\m x}$, naturally induces a monoid homomorphism $\mathscr{L}(V)\to \mathscr{L}(V_0)$, $\p \mapsto \overline \p$, which sends $\Image(\psi)$ into $\mathscr{L}^+(V_0)$ by assumption on the parameters $\m{a,\,b,\,c,\,d}$ and $\m f_{i,\a}$, $(i,\a) \in I\times \Phi^+$. If we compose again by the monoid homomorphism $\mathscr{L}^+(V_0) \to \Bin(\Phi^+)$, $\p \mapsto R_\p$, of remark \ref{des matrices positives aux relations}, we thus get a monoid homomorphism 
$$B^+ \to \Bin(\Phi^+), \ b \mapsto R_b = R_{\overline{\psi_b}}$$
\end{defi}

A quick look at the formulas of definition \ref{def psii} then easily gives the following : 

\begin{lem}[{\cite[4.3 and 4.4]{H}}]\label{Hee2}
Under condition (ii) of theorem \ref{criterion} and with the notations introduced in definition \ref{relations de Digne} above, the morphism $B^+ \to \Bin(\Phi^+)$, $b \mapsto R_b$, satisfies the following properties :
\begin{enumerate}
 \item $\a_i \not \in R_{\gras s_i}(\Phi^+)$,
\item if $i \neq j$, then $\a_i R_{\gras s_j} \a_i$,
\item if $m_{i,j} = 3$, then $\a_i R_{\gras s_j}R_{\gras s_i} \a_j$.
\end{enumerate}
\end{lem}

And a combination of the three previous lemmas easily gives theorem \ref{criterion}.

\subsection{Application of the faithfulness criterion}\label{Comments}\mbox{}\medskip

The typical situation where theorem \ref{criterion} applies is, for any totally ordered commutative ring $\m R_0$ (for example any subring of $\R$) : 
\begin{itemize}
\item $\m R = \m R_0[x]$ where $x$ is an indeterminate, and the evaluation at $x = 0$ for the morphism $\m R \to \m R_0$, 
\item $\m b,\, \m c,\, \m d$ positive units of $\m R_0$ with $\m{bc} < \m d^2$,
\item an LK-family $(f_i)_{i\in I}$ with $\Image(f_i) \subseteq x\m R$ and $f_i(e_{\a_i}) \neq 0$ for all $i \in I$.
\end{itemize}

Indeed, in that situation, condition (ii) of theorem \ref{criterion} is clearly satisfied. To see condition (i), note that, since the elements $f_i(e_{\a_i})$, $i\in I$, are non-zero, they become units of some appropriate overring $\m R'$ of $\m R$ (for example its field of fractions $\K(x)$, where $\K$ is the field of fractions of $\m R_0$). So if we denote by $V'$ the free $\m R'$-module with basis $(e_\a)_{\a \in \Phi^+}$, then $\Image(\psi)$ is included in $\EL(V) \cap \GL(V')$ and hence is cancellative. 

Moreover, the faithful LK-representation $\psi$ then induces an LK-representation $\psi_{\gr} : B \to \GL(V')$, which is faithful when $\G$ is spherical (in view of lemma \ref{inject monoid et group}). 

\medskip

Hence to construct faithful LK-representations of $B^+$ --- and of $B$ when $\G$ is spherical --- it suffices to construct LK-families over $x\m R_0[x]$ with non-zero elements $f_i(e_{\a_i})$ for $i\in I$. This will be done, for any Coxeter graph of small type with no triangle, and for the triangle graph $\tilde A_2$, in subsection \ref{Comments fin} below.

\begin{exa}Following \cite{K2, CW, Di, P}, one can choose $\m R_0 = \Z[y^{\pm 1}]$ for some $y \in \R_+^\star \setminus \{1\}$, and  $\m{(b,c,d)} = (y^p,y^q,y^r)$ with $p,\,q,\,r\in \Z$ such that $2r < p+q$ (resp. $2r > p+q$) if $0< y< 1$ (resp. $y > 1$).
\end{exa}

In \cite{K2, CW, Di, P}, the authors choose $0< y< 1$, $\m d = 1$ and $\m{(b,c)} = (y, y)$ (in \cite{CW}), $(1, y)$ (in \cite{Di}) or $(y, 1)$ (in \cite{P}). Note that the situation in \cite{K2} is slightly out of our settings since the value of $\m{(b,c)}$ varies for $\p_i$ (between $(y, 1)$ and $(1,y)$), depending on the considered $\{i\}$-mesh of cardinality two. The authors construct LK-families $(f_i)_{i\in I}$ over $x\Z[y]$ with elements $f_i(e_{\a_i})$, $i\in I$, all equal to $xy^4$ (in \cite{CW}) or $xy^2$ (in \cite{K2,Di,P}), hence the overring $\m R' = \Z[x^{\pm 1},y^{\pm 1}]$ is appropriate in the discussion above.

\section{On LK-families}\label{On LK-families}

We study in this section the module of LK-families $\F$ and its subset $\F_{\gr}$ for a given choice of parameters $\m{(b,c,d)} \in (\m R^\times)^3$, as defined in definition \ref{def lkp}. 

In subsection \ref{LK-families}, we characterize the LK-families in terms of relations between the elements $\m f_{i,\a} = f_i(e_\a) \in \m R$, for every $(i,\a) \in I\times \Phi^+$. This characterization generalizes some computations of \cite{CW, Di, P}, and we recall in subsections \ref{cas spherique} and \ref{The construction of Paris} below their results on those families : in our settings, it is shown in \cite{CW, Di} that, when $\G$ is spherical, there exists an isomorphism from $\F$ onto $\m R$ which sends $\F_{\gr}$ onto $\m R^\times$, an it is shown in \cite{P} that $\F_{\gr}$ is not trivial when $\G$ has no triangle.

The aim of subsection \ref{The affine case} is to explicit the structure of $\F$ and $\F_{\gr}$ when $\G$ is affine : we show that there exists an isomorphism from $\F$ onto $\m R^\N$ which sends $\F_{\gr}$ onto $\m R^\times\times \m R^{\N_{\geqslant 1}}$ (see theorem \ref{premier main thm}). In particular, this result holds for the affine type $\tilde A_2$ and hence gives the firsts examples of LK-representations of an Artin-Tits group whose type has triangles.

\medskip

In all this section, we fix a Coxeter matrix of small type $\G = (m_{i,j})_{i,j\in I}$, a commutative ring $\m R$, a triple $\m{(b,c,d)} \in (\m R^\times)^3$ and we set $\mathfrak a = \m d - \frac{\m{bc}}{\m d}$. We denote by $V$ the free $\m R$-module with basis $(e_\a)_{\a\in \Phi^+}$ and define the linear maps $\p_i \in \mathscr{L}(V)$, $i \in I$, as in definition \ref{def psii}. 

\subsection{Characterization of LK-families}\label{LK-families}\mbox{}\medskip

The following proposition gives a characterization of an LK-family $(f_i)_{i\in I}$ in terms of relations between the $\m f_{i,\a} = f_i(e_\a)$'s. This generalizes \cite[Prop. 3.2]{CW} and the computations of \cite[proof of Thm. 3.8]{Di} and \cite[proofs of lemmas 3.6 and 3.7]{P}.

\begin{pro}\label{relations for an LK-family}An element $(f_i)_{i\in I} \in (V^\star)^I$ is an LK-family if and only if the relations listed in \textsc{Table 1} below hold among the elements $\m f_{i,\a}$, $(i,\a) \in I\times \Phi^+$. (Note that the relations of cases \emph{(6)} and \emph{(8)} must hold whether $(\a_i|\a)$ is positive, zero or negative.)
\end{pro}

\begin{table}[ht]
 \begin{tabular}{|c|c|c|}\hline 
n$^\circ$ & Relations among the $\m f_{i,\a}$'s & Configuration of roots \\ \hline \hline
(1) & $\m f_{i,\a_j} = 0$ & if $i\neq j$ \\ \hline
(2) & $\m f_{i,\a_i} = \m f_{j,\a_j}$ & if $m_{i,j} = 3$ \\ \hline 
\begin{tabular}{c}(3)\\(4)\end{tabular} & \begin{tabular}{c}$\m f_{i,\b'} = \m f_{j,\b}$ \\ $\m{cf}_{i,\g'}+\m{af}_{i,\del} = \m{cf}_{j,\g}+\m{af}_{j,\del}$ \end{tabular}& if $m_{i,j} = 3$ and 
\begin{picture}(60,42)(-10,25)
\put(20,0){\circle*{4}}
\put(40,20){\circle*{4}}
\put(40,40){\circle*{4}}
\put(20,60){\circle*{4}}
\put(0,20){\circle*{4}}
\put(0,40){\circle*{4}}
\put(20,0){\line(1,1){20}}
\put(0,20){\line(0,1){20}}
\put(40,20){\line(0,1){20}}
\put(20,0){\line(-1,1){20}}
\put(20,60){\line(1,-1){20}}
\put(20,60){\line(-1,-1){20}}
\put(33,5){\small $j$}
\put(33,50){\small $j$}
\put(2,5){\small $i$}
\put(5,50){\small $i$}
\put(-6,26){\small $j$}
\put(42,26){\small $i$}
\put(-10,40){$\b$}
\put(44,40){$\b'$}
\put(-10,15){$\g$}
\put(44,15){$\g'$}
\put(23,-6){$\del$}
%\put(12,2){$\a$}
%\put(12,22){$\b$}
\end{picture} \\ & & \\ & & \\ \hline
(5) & $\m f_{i,\a} = \m f_{j,\a}$ & if $m_{i,j} = 3$ and 
\begin{picture}(35,20)(3,3)
\put(20,3){\circle*{4}}
\put(25,3){\circle{10}}
\put(32,0){\small $j$}
\put(15,3){\circle{10}}
\put(5,0){\small $i$}
\put(18,10){$\a$}
\end{picture} \\ & & \\ \hline \hline 
(6) & $\m{df}_{i,\a} = \m b\m f_{i,\b}$ & if $m_{i,j} = 2$ and 
\begin{picture}(35,20)(3,10)
%\put(20,0){\circle*{4}}
\put(20,0){\circle*{4}}
\put(20,20){\circle*{4}}
\put(20,0){\line(0,1){20}}
\put(22,8){\small $j$}
\put(12,0){$\b$}
\put(12,22){$\a$}
\end{picture} \\ & & \\  \hline
(7) & $\m c \m f_{i,\a_i+\a_j} = -\m a\m f_{i,\a_i}$ & if $m_{i,j} = 3$ \\ \hline
(8) & $\m c\m f_{i,\a} = \m b\m f_{j,\g} - \m a\m f_{i,\b}$ & if $m_{i,j} = 3$ and 
\begin{picture}(35,32)(3,20)
%\put(20,0){\circle*{4}}
\put(20,0){\circle*{4}}
\put(20,40){\circle*{4}}
\put(20,20){\circle*{4}}
\put(20,0){\line(0,1){40}}
\put(22,8){\small $i$}
\put(22,28){\small $j$}
\put(12,42){$\a$}
\put(12,22){$\b$}
\put(12,0){$\g$}
\end{picture} \\ & & \\ & & \\ \hline
\begin{tabular}{c}(9)\\(10)\end{tabular} & \begin{tabular}{c}$\m c\m f_{i,\b} = \m d\m f_{j,\g} - \m a\m f_{i,\g}$ \\ $\m d\m f_{i,\a} = \m b\m f_{j,\b}$ \end{tabular} & if $m_{i,j} = 3$ and 
\begin{picture}(35,32)(3,20)
%\put(20,0){\circle*{4}}
\put(20,0){\circle*{4}}
\put(20,40){\circle*{4}}
\put(20,20){\circle*{4}}
\put(20,0){\line(0,1){40}}
\put(22,8){\small $j$}
\put(22,26){\small $i$}
\put(12,42){$\a$}
\put(12,22){$\b$}
\put(12,0){$\g$}
\put(25,0){\circle{10}}
\put(31,-2){\small $i$}
\put(25,40){\circle{10}}
\put(31,38){\small $j$}
\end{picture} \\ & & \\ & & \\ \hline
\end{tabular} 

\vspace{1mm}

\textsc{Table 1}. Relations for an LK-family.
\end{table}

\proof This is simply the transcription on the basis elements $e_\a$, $\a\in \Phi^+$, of $V$, of conditions (i), (ii) and (iii) of definition \ref{def lkp} on the linear forms $f_i$, $i\in I$. Indeed, relation (1) is condition (i), and for every $i,\, j \in I$, we get by definition : 
\begin{center}
$\begin{cases}
 f_i\p_j(e_{\a}) = 0 & \text{  if  }\ \ \a = \a_j, \\
f_i\p_j(e_\a) = \m {df}_{i,\a} & \text{  if  
\begin{picture}(36,12)(0,-2)
\put(20,0){\circle*{4}}
\put(25,0){\circle{10}}
\put(32,-3){\footnotesize $j$}
\put(10,2){$\a$}
\end{picture}},\\
\begin{cases}f_i\p_j(e_\a) = \m b \m f_{i,\b} \\ f_i\p_j(e_\b) = \m a \m f_{i,\b} + \m c \m f_{i,\a}\end{cases} & \text{  if  
\begin{picture}(32,20)(0,8)
%\put(20,0){\circle*{4}}
\put(20,0){\circle*{4}}
\put(20,20){\circle*{4}}
\put(20,0){\line(0,1){20}}
\put(22,8){\footnotesize $j$}
\put(10,0){$\b$}
\put(10,18){$\a$}
\end{picture} in   } \, \Phi^+.
\end{cases}$
\end{center}
Assume that $m_{i,j} = 2$. We thus have $f_i\p_j(e_\a) = \m df_i(e_\a)$ if $\a = \a_j$ or if $\a$ is fixed by $s_j$. And for the $\{j\}$-mesh $\{\a,\b\}$ displayed above, 
$$\text{then  } \begin{cases}f_i\p_j(e_\a) = \m df_i(e_\a) \\ f_i\p_j(e_\b) = \m df_i(e_\b)\end{cases} \Longleftrightarrow \begin{cases}\m b \m f_{i,\b} = \m {d f}_{i,\a} \\ \m a \m f_{i,\b} + \m c \m f_{i,\a} = \m {d f}_{i,\b}\end{cases},$$ 
and both equations give relation (6) since $\dfrac{\m b}{\m d} = \dfrac{\m{d-a}}{\m c}$ in $\m R$.

Now assume that $m_{i,j} = 3$, and consider the system of equations $f_i\p_j(e_\a) = f_j\p_i(e_\a)$, for $\a$ running through the vertices of a given $\{i,j\}$-mesh $M$, and for the four possible types of $M$ (see remark \ref{mailles possibles}). 

Type 5 : $M = \{\a_i,\a_j,\a_i+\a_j\}$ (the situations for $\a_i$ and $\a_j$ are symmetrical), 
$$\text{then  } \begin{cases}f_i\p_j(e_{\a_i}) = f_j\p_i(e_{\a_i}) \\ f_i\p_j(e_{\a_i+\a_j}) = f_j\p_i(e_{\a_i+\a_j})\end{cases} \Longleftrightarrow \begin{cases}\m{af}_{i,\a_i}+ \m{cf}_{i,\a_i+\a_j} = 0 \\ \m{bf}_{i,\a_i} = \m{bf}_{j,\a_j}\end{cases},$$ this gives relations (7) and (2) (since $\m b \in \m R^\times$).

Type 6 : $M = \{\a\}$, $f_i\p_j(e_\a) = f_j\p_i(e_\a)\Longleftrightarrow \m{df}_{i,\a} = \m{df}_{j,\a}$, this is relation (5).

Type 7 : $M = \{\a,\b,\g\}$ as displayed above, 
\begin{center}
\begin{picture}(35,32)(3,15)
%\put(20,0){\circle*{4}}
\put(20,0){\circle*{4}}
\put(20,40){\circle*{4}}
\put(20,20){\circle*{4}}
\put(20,0){\line(0,1){40}}
\put(22,8){\small $j$}
\put(22,26){\small $i$}
\put(12,42){$\a$}
\put(12,22){$\b$}
\put(12,0){$\g$}
\put(25,0){\circle{10}}
\put(31,-2){\small $i$}
\put(25,40){\circle{10}}
\put(31,38){\small $j$}
\end{picture} then $\begin{cases}f_i\p_j(e_{\a}) = f_j\p_i(e_{\a}) \\ f_i\p_j(e_{\b}) = f_j\p_i(e_{\b}) \\ f_i\p_j(e_{\g}) = f_j\p_i(e_{\g}) \end{cases}\Longleftrightarrow \begin{cases}\m{df}_{i,\a}= \m{bf}_{j,\b} \\ \m{bf}_{i,\g}  = \m{af}_{j,\b} + \m{cf}_{j,\a}\\ \m{af}_{i,\g} + \m{cf}_{i,\b} = \m{df}_{j,\g}\end{cases},$
\end{center}
this gives relations (10), one case of (8) (by exchanging $i$ and $j$) and (9). 

Type 8 : $M = \{\a,\b, \b',\g,\g',\del\}$ as displayed above (the situations for $\b$ and $\b'$, and for $\g$ and $\g'$, are symmetrical),
\begin{center}
\begin{picture}(60,42)(-10,25)
\put(20,0){\circle*{4}}
\put(40,20){\circle*{4}}
\put(40,40){\circle*{4}}
\put(20,60){\circle*{4}}
\put(0,20){\circle*{4}}
\put(0,40){\circle*{4}}
\put(20,0){\line(1,1){20}}
\put(0,20){\line(0,1){20}}
\put(40,20){\line(0,1){20}}
\put(20,0){\line(-1,1){20}}
\put(20,60){\line(1,-1){20}}
\put(20,60){\line(-1,-1){20}}
\put(33,5){\small $j$}
\put(33,50){\small $j$}
\put(2,5){\small $i$}
\put(5,50){\small $i$}
\put(-6,26){\small $j$}
\put(42,26){\small $i$}
\put(-10,40){$\b$}
\put(44,40){$\b'$}
\put(-10,15){$\g$}
\put(44,15){$\g'$}
\put(23,-7){$\del$}
\put(22,63){$\a$}
%\put(12,2){$\a$}
%\put(12,22){$\b$}
\end{picture} then $\begin{cases}f_i\p_j(e_{\a}) = f_j\p_i(e_{\a}) \\ f_i\p_j(e_{\b}) = f_j\p_i(e_{\b}) \\ f_i\p_j(e_{\g}) = f_j\p_i(e_{\g}) \\ f_i\p_j(e_{\del}) = f_j\p_i(e_{\del}) \end{cases} \Longleftrightarrow \begin{cases}\m{bf}_{i,\b'}  = \m{bf}_{j,\b}  \\ \m{bf}_{i,\g} = \m{af}_{j,\b} + \m{cf}_{j,\a} \\ \m{af}_{i,\g} + \m{cf}_{i,\b} = \m{bf}_{j,\del} \\ \m{af}_{i,\del} + \m{cf}_{i,\g'}  = \m{af}_{j,\del} + \m{cf}_{j,\g} \end{cases},$
\end{center}
this gives relations (3) (since $\m b\in \m R^\times$), the two last cases of {(8)}, and {(4)}.\qed

\medskip 

Note that these relations are of two kinds : relations {(1)} to {(5)} give equalities between elements associated with roots of the same depth, whereas relations {(6)} to {(10)} express an element $\m f_{i,\a}$ in terms of a linear combination of some $\m f_{j,\b}$'s with $\dep(\b)< \dep(\a)$. In fact, relation {(5)} can be deleted from \textsc{Table 1} (this generalizes the analogous observation of \cite[proof of Prop. 3.2]{CW}, \cite[proof of Thm. 3.8]{Di}, and \cite[Lem. 3.5]{P}), this is proved in lemma \ref{pas de rel (5)} below.

\begin{lem}\label{cool}
 Let $i,\,j,\,k\in I$ be such that $m_{i,j} = m_{j,k} = m_{k,i} = 3$. Then for every $\a \in \Phi^+$, we have $(\a_i|\a) + (\a_j|\a) + (\a_k|\a) \leqslant 0$.
\end{lem}
\proof The map $v \mapsto (\a_i|v) + (\a_j|v) + (\a_k|v)$ is a linear form on $E = \oplus_{l\in I}\R\a_l$ which is clearly non-positive on the basis elements $\a_l$, $l\in I$. This gives the result since every $\a \in \Phi^+$ is a linear combination, with non-negative coefficients, of those elements $\a_l$, $l\in I$. \qed

\begin{lem}\label{pas de rel (5)}
Relation \emph{(5)} of \textsc{Table 1} is implied by relations \emph{(1)}, \emph{(6)}, \emph{(8)} and \emph{(10)}.
\end{lem}
\proof Fix $i,\, j \in I$ with $m_{i,j} = 3$ and $\a \in \Phi^+$ such that $(\a_i|\a)=(\a_j|\a) = 0$. Let us show by induction on $\dep(\a)$ that $\m f_{i,\a} = \m f_{j,\a}$, using only relations {(1)}, {(6)}, {(8)} and {(10)}. If $\dep(\a) = 1$, \ie if $\a = \a_k$ for some $k \in I$, then $k \neq i,\, j$ and the result is given by {(1)}. So assume that $\dep(\a) \geqslant 2$ and fix $k \in I$ such that $(\a_k|\a)>0$~; we set $\b = s_k(\a) \in \Phi^+$. 

If $m_{i,k} = m_{j,k} = 2$, then we get $(\a_i|\b)=(\a_j|\b) = 0$, whence $\m f_{i,\b} = \m f_{j,\b}$ by induction and hence $\m f_{i,\a} = \m f_{j,\a}$ by {(6)} since $\m d \in \m R^\times$. 

If $m_{i,k} = 2$ and $m_{j,k} = 3$, then the graph of $[\a]_{\{i,j,k\}}$ is the following : 

\begin{picture}(15,50)(0,26)
\put(0,0){\circle*{4}}
\put(0,20){\circle*{4}}
\put(0,40){\circle*{4}}
\put(0,60){\circle*{4}}
\put(0,0){\line(0,1){60}}
\put(-5,60){\circle{10}}
\put(5,60){\circle{10}}
\put(-5,40){\circle{10}}
\put(5,20){\circle{10}}
\put(-5,0){\circle{10}}
\put(5,0){\circle{10}}
\put(-7,58){\footnotesize $i$}
\put(3,58){\footnotesize $j$}
\put(-6,47){\footnotesize $k$}
\put(-7,38){\footnotesize $i$}
\put(-6,27){\footnotesize $j$}
\put(3,17){\footnotesize $k$}
\put(-6,7){\footnotesize $i$}
\put(-7,-2){\footnotesize $j$}
\put(3,-3){\footnotesize $k$}
\put(-3,66){$\a$}
\put(6,36){$\b$}
\put(-12,16){$\g$}
\put(-3,-12){$\delta$}
%\put(12,2){$\a$}
%\put(12,22){$\b$}
\end{picture}, whence 
\begin{tabular}{|lll}
 $\m f_{i,\a}\! \! \! \! $& $=\dfrac{\m b}{\m d}\m f_{i,\b}$ & by {(6)} \\
 & $=\dfrac{\overset{}{\m b}}{\m{cd}}\left({\m b}\m f_{j,\delta}-\m{a}\m f_{i,\g}\right)$ & by {(8)} \\
 $\m f_{j,\a}\! \! \! \! $ & $=\dfrac{\overset{}{1}}{\m c}\left(\m b\m f_{k,\g}-\m a\m f_{j,\b}\right)$ & by {(8)} \\
 & $=\dfrac{\overset{}{\m b}}{\m{cd}}\left({\m b}\m f_{k,\delta}-\m{a}\m f_{i,\g}\right)$ & by {(6)} and {(10)} \\ 
$\m f_{j,\delta}\! \! \! \! $ &  $=\overset{}{\m f_{k,\delta}}$ & by induction
\end{tabular}, and hence $\m f_{i,\a} = \m f_{j,\a}$.

\medskip 

Thanks to lemma \ref{cool}, we cannot have $m_{i,k} = m_{j,k} = 3$ in that situation, so we are done  (up to exchanging $i$ and $j$). \qed

\begin{rem}\label{rels irreductibilite}LK-families and reducibility.

If $\G_1,\ldots,\G_p$ are the connected components of $\G$, with vertex set $I_1,\ldots, I_p$ respectively, then $\Phi = \Phi_{I_1}\sqcup \cdots\sqcup \Phi_{I_p}$ and relations {(1)} and {(6)} imply that $\m f_{i,\a} = 0$ for every $(i,\a) \in I_m \times \Phi^+_{I_n}$ whenever $m \neq n$. As a consequence, any LK-representation $\psi$ of $B^+_\G$ is the direct sum of the induced LK-representations $\psi_n$ of $B^+_{\G_n}$ for $1\leqslant n \leqslant p$. 

Hence when considering LK-representations (or LK-families), there is no loss of generality in assuming that $\G$ is connected, in which case the elements $\m f_{i,\a_i}$, $i\in I$, are all equal by relation~{(2)}. %Recall however that in order to get invertible maps, we are interested in LK-families with $\m f_{i,\a_i} \in \m R^\times$ for $i \in I$. 
\end{rem}

\subsection{The spherical case}\label{cas spherique}\mbox{}\medskip

We assume here that $\G = A_n$ ($n \geqslant 1$), $D_n$ ($n\geqslant 4$), or $E_n$ ($n = 6,\, 7$ or $8$). In the following theorem, we rephrase the unicity statements of \cite[Prop. 3.5]{CW} and \cite[Thm. 3.8]{Di}. Recall that $\mathscr{F} \subseteq (V^\star)^I$ is the $\m R$-module of LK-families.

\begin{defi}Fix $i_0 \in I$. Let $\mu$ be the linear map $\mathscr{F} \to \m R, \ \, (f_i)_{i\in I} \mapsto \m f_{i_0,\a_{i_0}}$. In view of relation {(2)} (and of the fact that $\G$ is connected), $\mu$ does not depend on the choice of $i_0 \in I$.
\end{defi}

\begin{thm}The linear map $\mu : \mathscr{F} \to \m R$, is an isomorphism (or $\m R$-modules).
\end{thm}
\proof Since $\G$ is spherical, there is no mesh of type 8 (see remark \ref{mailles possibles}) in $\Phi^+$. 

Hence for a given LK-family $(f_i)_{i\in I}$, every $\m f_{i,\a}$ with $\dep(\a) \geqslant 2$ can be expressed as a linear combination of some $\m f_{j,\b}$'s with $\dep(\b) < \dep(\a)$, via at least one of the relations {(6)} to {(10)}. As a consequence, $(f_i)_{i\in I}$ is entirely determined by the values of the $\m f_{i,\a_j}$, for $i,\,j\in I$. And since $\m f_{i,\a_j} = 0$ if $i\neq j$ by {(1)}, and $\m f_{i,\a_i} = \m f_{i_0,\a_{i_0}}$ for every $i \in I$ by {(2)} (since $\G$ is connected), $(f_i)_{i\in I}$ is in fact entirely determined by the value $\m f_{i_0,\a_{i_0}}$, whence the injectivity of $\mu$. 

In order to show its surjectivity, the idea is to define an LK-family inductively, with basis step $\m f_{i,\a_i} = \m f\in \m R$ (one could chose $\m f = 1$ by linearity) and $\m f_{i,\a_j} = 0$ if $i\neq j$, and inductive step one of the suitable relations {(6)} to {(10)} to define $\m f_{i,\a}$ (with $\dep(\a)\geqslant 2$) in terms of a linear combination of some $\m f_{j,\b}$'s with $\dep(\b) < \dep(\a)$. Proving that the obtained family is indeed an LK-family amounts to proving that the definition of $\m f_{i,\a}$ does not depend on the choice of the suitable relation chosen in the inductive step. This is essentially done in \cite[Prop. 3.5]{CW} and \cite[Thm. 3.8]{Di}.\qed 

\medskip

When $\G$ is connected and spherical (and of small type), LK-representations of $B^+$ are then parametrized by $\m R$ and LK-representations of $B$ (those corresponding to LK-families with $\m f_{i_0,\a_{i_0}} \in \m R^\times$) are parametrized by $\m R^\times$. 

\medskip

Let us end this subsection with some consequences of that construction. 

\medskip

Note that since $\G$ is spherical, the free $\m R$-module $V$ is finite-dimensional and hence the notion of \emph{determinant} of an element of $\mathscr{L}(V)$ is defined.

\begin{cor}
Tow LK-representations $\psi$ and $\psi'$, associated with to distinct LK-families $(f_i)_{i\in I}$ and $(f'_i)_{i\in I}$ respectively, are non-equivalent.
\end{cor}
\proof It suffices to see that for a given $i \in I$, the maps $\psi_i$ and $\psi'_i$ have distinct determinant. But in view of remark \ref{les matrices}, we get $\det(\psi_i) = \m{uf}_{i,\a_i}$ and $\det(\psi'_i) = \m{uf}'_{i,\a_i}$ for a certain $\m u \in \m R^\times$, whence the result since the previous theorem shows that $(f_i)_{i\in I} \neq (f'_i)_{i\in I}$ implies $\m f_{i,\a_i} \neq \m f'_{i,\a_i}$. \qed

\medskip

Finally, an easy induction on $\dep(\a)$ gives the following remark, which generalizes \cite[Cor. 3.3]{CW} :

\begin{rem}Let $(f_i)_{i\in I}$ be an LK-family and set $f_{i_0,\a_{i_0}} = \m f$. 

Then for every $\a \neq \a_i$, we have $\m f_{i,\a} \in -\dfrac{\m{af}}{\m c}\m R$, and more precisely : 
\begin{enumerate}
\item $\m f_{i,\a} = 0$ if $i \not \in \Supp(\a)$, and 
\item $\m f_{i,\a} = -\dfrac{\m{af}}{\m c}\Big(\dfrac{\m b}{\m d}\Big)^{\dep(\a)-2}$ if $\dep(\a) \geqslant 2$ and $(\a_i|\a)>0$.
\end{enumerate}
 \end{rem}

\medskip

This construction can be generalized to an arbitrary Coxeter matrix of small type with no triangle, following \cite{P}. This is done in the following subsection.

\subsection{The LK-family of Paris}\label{The construction of Paris}\mbox{}\medskip

The main construction of \cite{P} is a uniform construction of an LK-family with $\m f_{i,\a_i} \in \m R^\times$ for every $i\in I$, for any Coxeter matrix of small type $\G = (m_{i,j})_{i,j\in I}$ with no triangle, \ie no subset $\{i,j,k\} \subseteq I$ with $m_{i,j} = m_{j,k} = m_{k,i} = 3$. 

\medskip 

This construction is made over $\m R = \Z[x^{\pm 1},y^{\pm 1}]$, with $\m f_{i,\a_i} = xy^2$ for every $i\in I$. The aim of this subsection is to generalize it to our settings. 
 
\begin{defi}[\cite{P}]\label{defpolParis}
 For every $\a \in \Phi^+$ with $\dep(\a) \geqslant 2$, fix an element $j_\a \in I$ such that $(\a|\a_{j_\a})>0$. Let us define a family $(\m f_{i,\a})_{(i,\a)\in I\times \Phi^+}$ by induction on $\dep(\a)$ as follows : 

\medskip

\noindent$\bullet$ Basis step : fix $\m f\in \m R$ (one can choose $\m f = 1$ by linearity) and set 

\begin{center}
\begin{tabular}{|c|c|c|}\hline Case & Value of $\m f_{i,\a}$ & Condition \\ \hline \hline 
$\overset{\mbox{}}{\text{(C1)}}$ & $\m f$ & if $\a = \a_i$ \\ \hline 
$\overset{\mbox{}}{\text{(C2)}}$ & $0$ & if $\a = \a_j$ for $j \neq i$ \\ \hline 
(C3) & $-\dfrac{\m{af}}{\m c}\Big(\dfrac{\m b}{\m d}\Big)^{\overset{}{\dep(\a)-2}}$ & if $\dep(\a) \geqslant 2$ and $(\a|\a_i)>0$ \\ \hline
\end{tabular}
\end{center}

\medskip

\noindent$\bullet$ Inductive step : if $\dep(\a) \geqslant 2$ and $(\a|\a_i)\leqslant 0$  --- hence $i \neq j_\a$ --- then set

\begin{center}
 \begin{tabular}{|c|c|c|}\hline Case & Value of $\m f_{i,\a}$ & Condition \\ \hline \hline 
(C4) & $\dfrac{\overset{}{\m b}}{\m d}\m f_{i,\b}$ & if $m_{i,j_\a} = 2$ \\ \hline & &  \\
(C5) & $\dfrac{1}{\m c}\left({\m b}\m f_{j_\a,\g}-{\m a}\m f_{i,\b}\right)$ & if $m_{i,j_\a} = 3$ and 
\begin{picture}(60,22)(-8,18)
\put(20,0){\circle*{4}}
\put(20,0){\circle*{4}}
\put(20,40){\circle*{4}}
\put(20,20){\circle*{4}}
\put(20,0){\line(0,1){40}}
\put(22,8){\footnotesize $i$}
\put(22,28){\footnotesize $j_\a$}
\put(12,42){$\a$}
\put(12,22){$\b$}
\put(12,0){$\g$}
\end{picture}
 \\ & & \\ & & \\ \hline & & \\
(C6) & $\dfrac{1}{\m c}\left({\m d}\m f_{j_\a,\b}-{\m a}\m f_{i,\b}\right)$ & if $m_{i,j_\a} = 3$ and 
\begin{picture}(60,22)(-8,18)
%\put(20,0){\circle*{4}}
\put(20,0){\circle*{4}}
\put(20,40){\circle*{4}}
\put(20,20){\circle*{4}}
\put(20,0){\line(0,1){40}}
\put(22,8){\footnotesize $j_\a$}
\put(22,26){\footnotesize $i$}
%\put(12,42){$\a$}
\put(12,22){$\a$}
\put(12,0){$\b$}
\put(25,0){\circle{10}}
\put(31,-2){\footnotesize $i$}
\put(25,40){\circle{10}}
\put(31,38){\footnotesize $j_\a$}
\end{picture} \\ & & \\ & & \\  \hline  
(C7) & $\dfrac{\m b}{\m d}\m{f}_{i,\b} + \dfrac{\m d}{\m c}\m f_{j_\a,\b} + \dfrac{\m{adf}}{\m c^2}\Big(\dfrac{\m b}{\m d}\Big)^{\dep(\a)-3}$ & if $m_{i,j} = 3$ and 
\begin{picture}(60,45)(-10,25)
\put(20,0){\circle*{4}}
\put(40,20){\circle*{4}}
\put(40,40){\circle*{4}}
\put(20,60){\circle*{4}}
\put(0,20){\circle*{4}}
\put(0,40){\circle*{4}}
\put(20,0){\line(1,1){20}}
\put(0,20){\line(0,1){20}}
\put(40,20){\line(0,1){20}}
\put(20,0){\line(-1,1){20}}
\put(20,60){\line(1,-1){20}}
\put(20,60){\line(-1,-1){20}}
\put(33,5){\footnotesize $j_\a$}
\put(33,50){\footnotesize $j_\a$}
\put(2,5){\footnotesize $i$}
\put(5,50){\footnotesize $i$}
\put(2,27){\footnotesize $j_\a$}
\put(42,26){\footnotesize $i$}
%\put(-10,40){$\a$}
%\put(44,40){$\a'$}
%\put(-10,15){$\b$}
\put(44,15){$\a$}
\put(23,-7){$\b$}
%\put(12,2){$\a$}
%\put(12,22){$\b$}
\end{picture} \\ & & \\ & & \\ & & \\ \hline
\end{tabular}
\end{center}

And we define the family $(f_i)_{i\in I} \in (V^\star)^I$ via $f_i(e_\a) = \m f_{i,\a}$ for $(i,\a) \in I\times \Phi^+$.
\end{defi}

Note that cases (C4) and (C5) occur whether $(\a_i|\a)$ is zero or negative. Case (C3) is a generalization of what happens when $\G$ is spherical, but is no longer a consequence of the relations of \textsc{Table 1} in general, and neither is case (C7) (see subsection \ref{The affine case}).

\begin{pro}
Assume that the family $(f_i)_{i\in I}$ of definition \ref{defpolParis} does not depend on the choice of the $j_\a$'s. Then it is an LK-family.
\end{pro}
\proof We have to show that the family $(\m f_{i,\a})_{(i,\a)\in I\times \Phi^+}$ satisfies the relations {(1)} to {(4)} and {(6)} to {(10)} of proposition \ref{relations for an LK-family} (thanks to lemma \ref{pas de rel (5)}). Relations {(1)}, {(2)}, {(3)}, {(7)} and {(10)} are clearly satisfied by construction. In the same way, (C3) implies relations {(6)} and {(8)} when $(\a_i|\a) > 0$ (use $\m{d(a-d) + bc} = 0$ to establish {(8)}).

Now consider a relation {(4)}, {(6)} with $(\a_i|\a) \leqslant 0$, {(8)} with $(\a_i|\a) \leqslant 0$, or {(9)}. Then the elements $\m f_{i,\a}$, for $\a$ of highest depth among the roots involved in this relation, are defined by induction. Under the assumption of the proposition, we are free to use the suitable case among (C7), (C4), (C5) or (C6) respectively, to define them ; this clearly shows that the considered relation is satisfied (use the fact that $\m{d(a-d) + bc} = 0$ to establish relation {(4)} via (C7)). \qed

\medskip

The previous proposition generalizes the computations of \cite[lemmas 3.5, 3.6 and 3.7]{P}. It is not clear whether the independance assumption is true in general, but this is at least the case when $\G$ has no triangle : 

\begin{pro}
 Assume that $\G$ has no triangle. Then the family $(f_i)_{i\in I}$ of definition \ref{defpolParis} does not depend on the choice of the $j_\a$'s.
\end{pro}
\proof This is \cite[lemmas 3.3 and 3.4]{P} : our settings are slightly more general, but the (long) computations of the proofs are essentially the same.\qed

\medskip

Hence if $\G$ has no triangle, then the module of LK-families $\mathscr{F}$ is not trivial (contains a free submodule of dimension 1), and, by choosing $\m f \in \m R^\times$ in definition \ref{defpolParis} above, one obtains an element of $\F_{\gr}$.

It can also be shown that the family of definition \ref{defpolParis} does not depend on the choice of the $j_\a$'s when $\G = \tilde A_2$ (following the same steps as in the proof of lemma \ref{cas pdelta moins alphai} below), hence the same holds for this triangle graph. We will more generally explicit all the LK-families for any affine Coxeter graph in the following section.

\subsection{The affine case}\label{The affine case}\mbox{}\medskip

We assume here that $\G = (m_{i,j})_{0 \leqslant i,j \leqslant n}$ is a Coxeter matrix of type $\tilde A_n$ ($n \geqslant 2$), $\tilde D_n$ ($n \geqslant 4$) or $\tilde E_n$ ($n = 6$, $7$, $8$). We set $I = [\![ 0,n]\!]$.

\medskip 

We denote by $\G_0 = (m_{i,j})_{1 \leqslant i,j \leqslant n}$ the corresponding spherical Coxeter matrix. Let $\Phi$ (resp. $\Phi_0$) be the root system associated with $\G$ (resp. $\G_0$) in $E = \oplus_{i=0}^n\R\a_i$ (resp. $E_0 = \oplus_{i=1}^n\R\a_i$) and let $\delta$ be the first positive imaginary root of $\Phi$, then we have the following decomposition (see~\cite{Kac}) : 
$$\Phi = \bigsqcup_{p \in \Z}\big(\Phi_{0}+p\delta\big)\ \ \text{  and  } \ \ \Phi^+ = \Phi_0^+ \bigsqcup \Big(\bigsqcup_{p \in \N_{\geqslant 1}}\big(\Phi_{0}+p\delta\big)\Big).$$

As a consequence, we get the following remark :

\begin{rem}\label{hexagones affines}The only meshes of type 8 in $\Phi^+$ (see remark \ref{mailles possibles}) are the following ones, for $p\geqslant 1$ and $m_{i,j} = 3$ : 
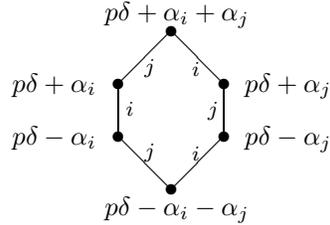
\begin{center}
 \begin{picture}(100,60)(-30,5)
\put(20,0){\circle*{4}}
\put(0,20){\circle*{4}}
\put(40,20){\circle*{4}}
\put(0,40){\circle*{4}}
\put(40,40){\circle*{4}}
\put(20,60){\circle*{4}}
\put(20,0){\line(1,1){20}}
\put(20,0){\line(-1,1){20}}
\put(0,20){\line(0,1){20}}
\put(40,20){\line(0,1){20}}
\put(20,60){\line(1,-1){20}}
\put(20,60){\line(-1,-1){20}}
\put(28,11){\footnotesize $i$}
\put(28,43){\footnotesize $i$}
\put(34,28){\footnotesize $j$}
\put(3,28){\footnotesize $i$}
\put(10,12){\footnotesize $j$}
\put(10,44){\footnotesize $j$}
\put(-5,-10){$p\delta-\a_i-\a_j$}
\put(48,17){$p\del-\a_j$}
\put(48,37){$p\del+\a_j$}
\put(-40,17){$p\del-\a_i$}
\put(-40,37){$p\del+\a_i$}
\put(-5,64){$p\del+\a_i+\a_j$}
\end{picture}
\end{center}\vspace{5mm}
In particular, for a given $(i,\a) \in I\times\Phi^+$ with $\dep(\a) \geqslant 2$, then either $\a = p\del \pm \a_i$ for some $p \geqslant 1$, or the pair $(i,\a)$ appears at the left-hand side of (at least) one of the relations {(6)} to {(10)} of \textsc{Table 1}, and for every such relation and every pair $(j,\b)$ involved in its right-hand side, then $\b \neq q\del - \a_j$ for every $q \geqslant 1$.
\end{rem}

\begin{defi}\label{def mu affine}Let $i_0,\, j_0 \in I$ be such that $m_{i_0,j_0} = 3$. We denote by $\mu$ the linear map 
$\mathscr{F} \to \m R^\N, \ \, (f_i)_{i\in I} \mapsto (\m f_n)_{n\in \N}$, where 
\begin{enumerate}
\item $\m f_{2p} = \m f_{i_0,p\del+\a_{i_0}}$ for every $p\in \N$, 
\item $\m f_{2p-1} = \m{cdf}_{i_0,p\del-\a_{i_0}} - \m{bcf}_{i_0,p\del-\a_{i_0}-\a_{j_0}} - \m{d}^2\m f_{j_0,p\del-\a_{i_0}-\a_{j_0}}$ for every $p\in \N_{\geqslant 1}$.
\end{enumerate}
\end{defi}

We will show in proposition \ref{mu indep de i0 j0} below that $\mu$ does not depend on the choice of $i_0,\, j_0 \in I$ such that $m_{i_0,j_0} = 3$. The aim of this subsection is then to prove the following : 

\begin{thm}\label{premier main thm}The linear map $\mu : \mathscr{F} \to \m R^\N$ is an isomorphism (of $\m R$-modules).
\end{thm}

Hence in those cases, LK-representations of $B^+$ are parametrized by $\m R^\N$ and LK-representations of $B$ (those corresponding to LK-families with $\m f_{i_0,\a_{i_0}} \in \m R^\times$) are parametrized by $\m R^\times\times \m R^{\N_{\geqslant 1}}$. The injectivity and surjectivity of $\mu$ are proved respectively in propositions \ref{condition necessaire} and \ref{condition suffisante} below. 

\begin{nota}For every $k\in \N_{\geqslant 1}$, let us denote by $\Phi^+_{k}$ the subset of $\Phi^+$ composed of the positive roots of depth smaller than (or equal to) $k$. 
\end{nota}

In the two following lemmas, we assume that we are given a family $\m F_k = (\m f_{i,\a})_{(i,\a) \in I \times \Phi^+_{k}} \in \m R^{I\times \Phi^+_{k}}$ whose elements satisfy the relations of \textsc{Table 1} whenever the roots involved are of depth smaller than (or equal to) $k$, and it is understood that we work with the elements of $\m F_k$.

\begin{lem}\label{cas pdelta moins alphai}Fix $(i,p) \in I\times\N_{\geqslant 1}$ and assume that $k = \dep(p\del -\a_i)-1$. Then the element $\m{bcf}_{i,p\del-\a_i-\a_j}+\m{d}^2\m f_{j,p\del - \a_i-\a_j}$
does not depend on $j\in I$ such that $m_{i,j} = 3$.
\end{lem}
\proof Assume that $j,\,k \in I$ are such that $m_{i,j} = m_{i,k} = 3$. 

If $m_{j,k} = 2$, then the result follows from relations {(6)} and {(9)} : indeed, we get 
 \begin{center}
\begin{tabular}{rcll}
$\m{cf}_{i,p\del-\a_i-\a_j}$ & $=$ & $\m{df}_{k,p\del-\a_i-\a_j-\a_k} - \m{af}_{i,p\del-\a_i-\a_j-\a_k}$ & by {(9)}, \\
$\m{cf}_{i,p\del-\a_i-\a_k}$ & $=$ & $\m{df}_{j,p\del-\a_i-\a_j-\a_k} - \m{af}_{i,p\del-\a_i-\a_j-\a_k}$ & by {(9)}, \\
$\m{df}_{j,p\del - \a_i-\a_j}$ & $=$ & $\m{bf}_{j,p\del - \a_i-\a_j-\a_k}$ & by {(6)},\\
$\m{df}_{k,p\del - \a_i-\a_k}$ & $=$ & $\m{bf}_{k,p\del - \a_i-\a_j-\a_k}$ & by {(6)}.
\end{tabular}
\end{center}

If $m_{j,k} = 3$, then $\G = \tilde A_2$, $\{i,j,k\} = \{0,1,2\}$ and $p\del-\a_i-\a_j = (p-1)\del+\a_k$. In that case we can prove more, namely that the value of $\m f_{l,(p-1)\del+\a_m}$ does not depend on the pair $(l,m)\in \{0,1,2\}^2$ such that $l \neq m$. To do this, one can first prove the similar statement for the elements $\m f_{l,q\del-\a_m}$ with $1 \leqslant q \leqslant p-1$ by induction on $q$, thanks to relations {(3)} and {(8)} (the case $q = 1$ is given by relations {(7)} and {(2)}), and then prove the desired statement for the elements $\m f_{l,q\del+\a_m}$ with $0 \leqslant q \leqslant p-1$ by induction on $q$, thanks to relation {(8)} and the intermediate result (the case $q = 0$ is given by relation {(1)}). \qed 

\begin{lem}\label{cas pas pdelta moins alphai}Fix $(i,\a)\in I\times\Phi^+$, with $\a \neq p\del \pm \a_i$ for every $p\in \N$, and assume that $k = \dep(\a)-1$. If we define $\m f_{i,\a} \in \m R$ by one of the relations \emph{(6)} to \emph{(10)} where $(i,\a)$ appears at the left-hand side, then the value of $\m f_{i,\a}$ does not depend on the chosen relation.
\end{lem}
\proof Let us assume that $(i,\a)$ appears at the left-hand side of two of the relations {(6)} to {(10)} and let us denote by $j$ (resp. $k$) the index distinct from $i$ involved in the first (resp. the second) of those two relations. Note that, in that situation, then $\G \neq \tilde A_2$ and there are only ten possibilities up to exchanging $j$ and $k$ (use remark \ref{hexagones affines}) : two relations {(6)} with $m_{j,k} = 2$ or $m_{j,k} = 3$, relations {(6)} and {(8)} with $m_{j,k} = 2$, relations {(6)} and {(9)} with $m_{j,k} = 2$ or $m_{j,k} = 3$, relations {(6)} and {(10)} with $m_{j,k} = 2$ or $m_{j,k} = 3$, two relations {(8)} with $m_{j,k} = 2$, two relations {(9)} with $m_{j,k} = 2$, two relations {(10)} with $m_{j,k} = 2$. 

For example in the case {(6)} and {(9)} with $m_{j,k} = 2$, the graph of $[\a]_{\{i,j,k\}}$ is 
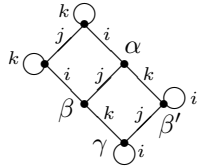
\begin{center}
\begin{picture}(50,55)(0,-5)
\put(0,30){\circle*{3}}
\put(15,15){\circle*{3}}
\put(15,45){\circle*{3}}
\put(30,0){\circle*{3}}
\put(30,30){\circle*{3}}
\put(45,15){\circle*{3}}
\put(30,0){\line(-1,1){30}}
\put(45,15){\line(-1,1){30}}
\put(30,0){\line(1,1){15}}
\put(15,15){\line(1,1){15}}
\put(0,30){\line(1,1){15}}
\put(-4,30){\circle{8}}
\put(15,49){\circle{8}}
\put(49,17){\circle{8}}
\put(30,-4){\circle{8}}
\put(7,24){\scriptsize $i$}
\put(4,39){\scriptsize $j$}
\put(19,24){\scriptsize $j$}
\put(22,39){\scriptsize $i$}
\put(22,9){\scriptsize $k$}
\put(34,9){\scriptsize $j$}
\put(37,24){\scriptsize $k$}
\put(-14,30){\scriptsize $k$}
\put(5,48){\scriptsize $k$}
\put(55,15){\scriptsize $i$}
\put(35,-5){\scriptsize $i$}
\put(30,34){$\a$}
\put(5,10){$\b$}
\put(18,-3){$\g$}
\put(43,4){$\b'$}
\end{picture}
\end{center}
whence $\left\{\begin{array}{rcll} \dfrac{\m b}{\m d}\m f_{i,\b} & = &\dfrac{\m b}{\m c}\m f_{k,\g}  - \dfrac{\m {ab}}{\m{cd}}\m f_{i,\g}& \text{by {(9)}, and} \\ \dfrac{\m {d}}{\m{c}}\m f_{k,\b'} - \dfrac{\m {a}}{\m{c}}\m f_{i,\b'} & = & \dfrac{\overset{}{\m b}}{\m c}\m f_{k,\g} - \dfrac{\m {ab}}{\m{cd}}\m f_{i,\g}& \text{by {(6)} (two times)}. \end{array}\right.$

And in the case {(6)} and {(9)} with $m_{j,k} = 3$, the graph of $[\a]_{\{i,j,k\}}$ is 
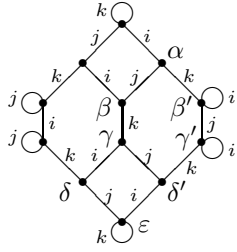
\begin{center}
\begin{picture}(65,80)
 \put(0,30){\circle*{3}}
 \put(0,45){\circle*{3}}
 \put(15,15){\circle*{3}}
 \put(15,60){\circle*{3}}
 \put(30,0){\circle*{3}}
 \put(30,30){\circle*{3}}
 \put(30,45){\circle*{3}}
 \put(30,75){\circle*{3}}
 \put(45,15){\circle*{3}}
 \put(45,60){\circle*{3}}
 \put(60,30){\circle*{3}}
 \put(60,45){\circle*{3}}
 \put(30,0){\line(-1,1){30}}
 \put(30,0){\line(1,1){30}}
 \put(30,75){\line(-1,-1){30}}
 \put(30,75){\line(1,-1){30}}
 \put(0,30){\line(0,1){15}}
 \put(30,30){\line(0,1){15}}
 \put(60,30){\line(0,1){15}}
 \put(30,30){\line(1,-1){15}}
 \put(30,30){\line(-1,-1){15}}
 \put(30,45){\line(1,1){15}}
 \put(30,45){\line(-1,1){15}}
 \put(30,-4){\circle{8}}
 \put(30,79){\circle{8}}
 \put(64,28){\circle{8}}
 \put(64,47){\circle{8}}
 \put(-4,30){\circle{8}}
 \put(-4,45){\circle{8}}
 \put(2,35){\scriptsize $i$}
 \put(32,35){\scriptsize $k$}
 \put(62,35){\scriptsize $j$}
 \put(-13,30){\scriptsize $j$}
 \put(-13,45){\scriptsize $j$}
 \put(70,25){\scriptsize $i$}
 \put(70,45){\scriptsize $i$}
 \put(20,78){\scriptsize $k$}
 \put(20,-8){\scriptsize $k$}
 \put(38,68){\scriptsize $i$}
 \put(18,68){\scriptsize $j$}
 \put(3,53){\scriptsize $k$}
 \put(53,53){\scriptsize $k$}
 \put(33,53){\scriptsize $j$}
 \put(23,53){\scriptsize $i$}
 \put(38,23){\scriptsize $j$}
 \put(18,23){\scriptsize $i$}
 \put(8,23){\scriptsize $k$}
 \put(54,18){\scriptsize $k$}
 \put(33,8){\scriptsize $i$}
 \put(23,8){\scriptsize $j$}
 \put(47,62){$\a$}
 \put(20,40){$\b$}
 \put(48,40){$\b'$}
 \put(21,30){$\g$}
 \put(50,29){$\g'$}
 \put(6,8){$\delta$}
 \put(47,8){$\delta'$}
 \put(36,-4){$\varepsilon$}
\end{picture}
\end{center}
\medskip
whence $\left\{\begin{array}{rcll} \dfrac{\m b}{\m d}\m f_{i,\b} & = &\dfrac{\m b^2}{\m{cd}}\m f_{k,\delta} - \dfrac{\m{ab}}{\m{cd}}\m f_{i,\g}& \text{by {(8)},}\\ & = & \dfrac{\m{b}^2}{\m{c}^2}\m f_{j,\varepsilon} -\dfrac{\m{ab}^2}{\m{c}^2\m d}\m f_{k,\varepsilon} - \dfrac{\m{ab}^2}{\m{cd}^2}\m f_{i,\delta'} & \text{by {(9)} and {(6)},}\\ \multicolumn{4}{l}{\dfrac{\m {d}}{\m{c}}\m f_{k,\b'} - \dfrac{\m {a}}{\m{c}}\m f_{i,\b'}} \\ & = & \dfrac{\m{bd}}{\m{c}^2}\m f_{j,\delta'} - \dfrac{\m{ad}}{\m{c}^2}\m f_{k,\g'} - \dfrac{\m{ab}}{\m{cd}}\m f_{i,\g'} & \text{by {(8)} and {(6)},}\\ & = & \dfrac{\m{b}^2}{\m{c}^2}\m f_{j,\varepsilon} - \dfrac{\m{ab}}{\m{c}^2}\m f_{i,\delta'} -\dfrac{\m{ab}^2}{\m{c}^2\m d}\m f_{k,\varepsilon} + \dfrac{\m{a}^2\m b}{\m{c}^2\m d}\m f_{i,\delta'}& \text{by {(6)}, {(10)} and {(8)},} \\ & = & \dfrac{\m{b}^2}{\m{c}^2}\m f_{j,\varepsilon} - \dfrac{\m{ab}^2}{\m{c}^2\m d}\m f_{k,\varepsilon} + \dfrac{\m{ab(a-d)}}{\m{c}^2\m d}\m f_{i,\delta'}, & \end{array}\right.$ 

\noindent and the result since $\frac{\m{ab(a-d)}}{\m{c}^2\m d} +\frac{\m{ab}^2}{\m{cd}^2} = \frac{\m{ab}}{\m{(cd)}^2}\m{(d(a-d)+bc)} = 0$.

The eight remaining cases are similar and left to the reader. \qed

\medskip

We are now able to prove the announced results on the linear map $\mu : \mathscr{F} \to \m R^\N$ of definition \ref{def mu affine}.

\begin{pro}\label{mu indep de i0 j0}
The definition of $\mu$ does not depend on the choice of $i_0,\, j_0 \in I$ such that $m_{i_0,j_0} = 3$.
\end{pro}
\proof Let $(f_i)_{i\in I}$ be an LK-family. Relations {(2)} and the fact $\G$ is connected show that the elements $\m f_{i,\a_i}$, $i\in I$, are all equal to $\m f_{0}$. Now fix $p \in \N_{\geqslant 1}$. In the same vein, relation {(3)} and the fact $\G$ is connected show that the elements $\m f_{i,p\del+\a_i}$, $i\in I$, are all equal to $\m f_{2p}$. 

Now if we set $\m f_{i,j} = \m{cdf}_{i,p\del-\a_{i}} - \m{bcf}_{i,p\del-\a_{i}-\a_{j}} - \m{d}^2\m f_{j,p\del-\a_{i}-\a_{j}}$ for $i,\, j\in I$ with $m_{i,j} = 3$, we are left to show that $\m f_{i,j} = \m f_{i_0,j_0} = \m f_{2p-1}$ for every $i,\, j\in I$ with $m_{i,j} = 3$. But lemma \ref{cas pdelta moins alphai} gives $\m f_{i,j} = \m f_{i,k}$ if $m_{i,j} = m_{j,k} = 3$, and since $\m a = \m d - \frac{\m{bc}}{\m d}$, relation {(4)} can be written $\frac{1}{\m d}\m f_{i,j} = \frac{1}{\m d}\m f_{j,i}$, whence $\m f_{i,j} = \m f_{j,i}$, when $m_{i,j} = 3$. The connectedness of $\G$ then gives the result. \qed

\begin{pro}\label{condition necessaire}The linear map $\mu$ is injective.
\end{pro}
\proof Let $(f_i)_{i\in I}$ be an LK-family. Thanks to remark \ref{hexagones affines}, an easy induction on $\dep(\a)$ shows that, for every $(i,\a)\in[\![ 0,n]\!]\times\Phi^+$, either $\a = p\del -\a_i$ for some $p \in \N_{\geqslant 1}$, or $\m f_{i,\a}$ is a linear combination of some $\m f_{j,p\del + \a_j}$'s for $j\in I$ and $p\in \N$. Hence $(f_i)_{i\in I}$ is entirely determined by the elements $\m f_{i,p\del \pm \a_i}$, for $i \in I$ and $p \in \N$. 

But by proposition \ref{mu indep de i0 j0}, the elements $\m f_{i,p\del + \a_i}$, $i \in I$, are all equal to $\m f_{2p}$ and we get $\m{f}_{i,p\del-\a_{i}} = \frac{\m b}{\m d}\m{f}_{i,p\del-\a_{i}-\a_{j}} + \frac{\m d}{\m c}\m f_{j,p\del-\a_{i}-\a_{j}} + \frac{1}{\m{cd}}\m f_{2p-1}$ for any $j \in I$ with $m_{i,j} = 3$, hence $\m{f}_{i,p\del-\a_{i}}$ is entirely determined by some $\m f_{2q}$, $q\in \N$, and $\m f_{2p-1}$. \qed

% \medskip
% 
% The family $(T_q)_{q \in \N}$ can in fact be chosen arbitrarily (with $T_0 \neq 0$) :
% 
\begin{pro}\label{condition suffisante}The linear map $\mu$ is surjective.
\end{pro}
\proof Fix $(\m f_n)_{n\in \N} \in \m R^\N$, and let us construct a family $(f_i)_{i\in I} \in (V^\star)^I$ by induction as follows (recall that we set $f_i(e_\a) = \m f_{i,\a}$ for every $(i,\a) \in I\times\Phi^+$) : 
\begin{itemize}
\item Basis step : we set $\m f_{i,\a_j} = 0$ for $i\neq j$, and $\m f_{i,p\del+\a_i} = \m f_{2p}$ for every $i\in I$ and $p\in \N$.
\item Inductive step : if $(i,\a)\in I\times\Phi^+$ is not handled by the basis step (hence $\dep(\a) \geqslant 2$) and is such that all the $\m f_{j,\b}$, for $j \in I$ and $\dep(\b) < \dep(\a)$, are constructed. Then 
\begin{enumerate}
\item if $\a = p\del-\a_i$, we set $\m f_{i,\a} = \dfrac{\m b}{\m d}\m{f}_{i,p\del-\a_{i}-\a_{j}} + \dfrac{\m d}{\m c}\m f_{j,p\del-\a_{i}-\a_{j}} + \dfrac{1}{\m{cd}}\m f_{2p-1}$ for some $j$ such that $m_{i,j} = 3$,
\item if not, then $(i,\a)$ appears at the left-hand side of (at least) one of the relations {(6)} to {(10)} (see remark \ref{hexagones affines}) ; we define $\m f_{i,\a}$ via the corresponding right-hand side.
\end{enumerate}
\end{itemize}

We are left to show that $(f_i)_{i\in I} \in (V^\star)^I$ is an LK-family, since it will then be, by construction, an antecedent of $(\m f_n)_{n\in \N} \in \m R^\N$ by $\mu$. We proceed by induction on $m \in \N$ in order to show that the relations of \textsc{Table 1} that involve only roots of depth smaller than (or equal to) $m$ are satisfied by the elements $\m f_{i,\a}$, for $i \in I$ and $\a \in \Phi^+$ with $\dep(\a) \leqslant m$. 

If $m = 0$, the only relations to consider are relations {(1)} and {(2)}, which are satisfied by construction of the basis step. Relation {(3)} is also satisfied for arbitrary depths by construction of the basis step. Now assume that we know the result for some $m \in \N$ and consider a relation that involves a root of depth $m+1$ and no root of higher depth. 

Assume first that it is a relation of type {(4)}, involving the indices $i$ and $j$ (and hence the roots $p\del - \a_i$, $p\del - \a_j$ and $p\del - \a_i- \a_j$). Lemma \ref{cas pdelta moins alphai} shows that the definition of $\m f_{i,p\del - \a_i}$ (resp. $\m f_{j,p\del - \a_j}$) at inductive step (i) does not depend on the choice of $k$ (resp. $k'$) such that $m_{i,k} = 3$ (resp. $m_{j,k'} = 3$). So we are free to chose $k = j$ (resp. $k' = i$), and we obtain that both sides of relation {(4)} are equal to $\m d\left(\m f_{i,p\del -\a_i - \a_j} + \m f_{j,p\del -\a_i - \a_j}\right) + \frac{1}{\m d}\m f_{2p-1}$ since $\mathfrak {a + \frac{bc}{d} = d}$.

Assume finally that it is a relation of type {(6)--(10)}. Lemma \ref{cas pas pdelta moins alphai} shows that the definition at inductive step (ii) does not depend on the choice of the relation, hence we are free to use the considered relation at this step and this gives the result. \qed 

\begin{exa}The LK-family of Paris (see definition \ref{defpolParis}) is the one corresponding to the family $(\m f_{n})_{n\in \N}$ with $\m f_0 = \m f$ and, for $p \geqslant 1$, $\m f_{2p} = -\dfrac{\m{af}}{\m c}\Big(\dfrac{\m b}{\m d}\Big)^{\dep(p\del + \a_i)-2}$ and $\m f_{2p-1} = \dfrac{\m{ad}^2\m{f}}{\m c}\Big(\dfrac{\m b}{\m d}\Big)^{\dep(p\del -\a_i)-3}$.
\end{exa}

\begin{exa}\label{cas Antilde}Let us assume that $\G = \tilde A_n$. Then each $\a \in \Phi^+$ has a unique decomposition as $\a = p\del + \sum_{k=j}^{j+\ell}\a_{\overline k}$, with $p \in \N$, $j \in [\![ 0,n]\!]$, $\ell \in [\![ 0,n-1]\!]$, and $\overline k$ the rest of $k$ modulo $n+1$. We then call {\it domain} (resp. {\it interior}, resp. {\it boundary}) of $\a$ the set $\overline \a = \{\overline k \mid j\leqslant k \leqslant j+\ell\}$ (resp. $\a^\circ = \{\overline k \mid j+1\leqslant k \leqslant j+\ell-1\}$, resp. $\partial\a = \overline \a \setminus \a^\circ = \{\overline j,\overline{j+\ell}\}$). If $(f_i)_{i\in I}$ is an LK-family, then one can check, by induction on $\dep(\a)$, that the element $\m f_{i,\a} = f_i(e_\a)$ is equal to : 
\begin{itemize}
\item $\m f_{2p}$ if $\a = p\del + \a_i$ (\ie $i \in \partial \a$ and $\ell = 0$), %or\\
\item $-\dfrac{\m a}{\m c}\Big(\dfrac{\m b}{\m d}\Big)^{\ell-1}{\displaystyle \sum_{q = 0}^p}\Big(\dfrac{\m b^{n-1}}{\m{cd}^{n-2}}\Big)^{p-q}\m f_{2q}$ if $i \in \partial \a$ and $\ell\geqslant 1$, 
\item $\Big(\dfrac{\m a}{\m c}\Big)^2\Big(\dfrac{\m b}{\m d}\Big)^{\ell-2}{\displaystyle \sum_{q = 0}^p}(p-q+1)\Big(\dfrac{\m b^{n-1}}{\m{cd}^{n-2}}\Big)^{p-q}\m f_{2q}$ if $i \in \a^\circ$,
\item $\dfrac{\m a^2}{\m{c}}\Big(\dfrac{\m b}{\m d}\Big)^{\ell}{\displaystyle \sum_{q = 0}^{p-1}}(p-q)\Big(\dfrac{\m b^{n-1}}{\m{cd}^{n-2}}\Big)^{p-q}\m f_{2q}$ if $i \not \in \overline \a$ and $\ell \leqslant n-2$, %or\\
\item $\displaystyle \dfrac{\m a^2}{\m{c}}\Big(\dfrac{\m b}{\m{d}}+\dfrac{\m d}{\m{c}}\Big)\Big(\dfrac{\m b}{\m d}\Big)^{n-2}\sum_{q = 0}^{p-1}(p-q)\Big(\dfrac{\m b^{n-1}}{\m{cd}^{n-2}}\Big)^{p-q}\m f_{2q} + \dfrac{1}{\m{cd}}\m f_{2p-1}$ if $\a = p\del -\a_i$ \\ (\ie $i \not \in \overline \a$ and $\ell = n-1$).
\end{itemize}
\end{exa}

\subsection{Comments}\label{Comments fin}\mbox{}\medskip

Assume that we are in the situation of subsection \ref{Comments}, that is, in particular, with $\m R = \m R_0[x]$ for some totally ordered commutative ring $\m R_0$. 

Then by choosing for $\m f$ an element of $x\m R$ for the basis step of the inductive construction of LK-families in subsections \ref{cas spherique} and \ref{The construction of Paris}, it is clear that the obtained LK-family $(f_i)_{i\in i}$ is such that $\Image(f_i) \subseteq x\m R$ for all $i \in I$. We get the same result in the affine case by choosing for $(\m f_n)_{n\in \N}$ a family of elements of $x\m R$ in the inductive construction of subsection \ref{The affine case}. 

If moreover $\m f$ (resp. $\m f_0$) is chosen to be non-zero, hence is a unit of some overring $\m R'$ of $\m R$, then the obtained LK-family will be suitable to apply the faithfulness criterion (theorem \ref{criterion}) to the associated LK-representation $\psi$ of $B^+$ in order to show that it is faithful. 

\section{Twisted LK-representations}\label{Action of graph automorphisms}

In \cite{Di}, Digne defines ``twisted'' LK-representations for an Artin-Tits group of non-small crystallographic and spherical type (\ie of type $B_n$, $F_4$ or $G_2$), using the fact that this group is the subgroup of fixed elements of an Artin-Tits group of small and spherical type ($A_{2n-1}$, $E_6$ or $D_4$ respectively) under a graph automorphism, and shows that those representations are faithful.

The aim of this section is to generalize this construction and the faithfulness result to any Artin-Tits monoid that appears as the submonoid of fixed elements of an Artin-Tits monoid of small type under a group of graph automorphisms. Note that our proof of faithfulness (\cf subsection \ref{Faithfulness of the induced representation}) is different from the one of \cite{Di} as it is general and avoid any case-by-case analysis.

\medskip

Let $\G = (m_{i,j})_{i,j\in I}$ be a Coxeter matrix of small type and let $G$ be a subgroup of $\Aut(\G)$. We fix $\m{(b,c,d)} \in (\m R^\times)^3$ and set $\m a = \m d -\frac{\m{bc}}{\m d}$. Let us consider an LK-family $(f_i)_{i\in I} \in \F$. Then we get an LK-representation $\psi : B^+ \to \mathscr{L}(V)$, $\gras s_i \mapsto \psi_i = \p_i + f_i\bt e_{\a_i}$, which induces an LK-representation of $B$ (\ie which has invertible images) whenever $(f_i)_{i\in I} \in \F_{\gr}$, \ie $f_i(e_{\a_i}) \in \m R^\times$ for every $i \in I$.

\subsection{Definition}\label{def twisted LK-reps}\mbox{}\medskip

Recall that the group $G$ naturally acts on $B^+$ and on $\Phi^+$ (see section \ref{Graph automorphisms}). The action of $G$ on $\Phi^+$ induces an action of $G$ on $V$ by permutation of the basis $(e_\a)_{\a \in \Phi^+}$. We denote by $(B^+)^G$ (resp. $V^G$) the submonoid (resp. submodule) of fixed points of $B^+$ (resp. of $V$) under the action of $G$. Recall that $(B^+)^G = B^+_{\G'}$ for a certain Coxeter graph $\G'$.

\medskip

We denote by $\p : B^+ \to \mathscr{L}(V)$, $b \mapsto \p_b$, the LK-representation of $B^+$ associated with the trivial LK-family (\ie with $\psi_i = \p_i$ for all $i\in I$).

\begin{lem}\label{stabilisation phi}For all $(b,v,g) \in B^+\times V \times \Aut(\G)$, we get $g\left(\p_b(v)\right) = \p_{g(b)}(g(v))$. In particular, for all $b \in (B^+)^G$, $\p_b$ stabilizes $V^G$ and hence $\p$ induces a linear representation $\p^G : (B^+)^G \to \mathscr{L}(V^G)$, $b \mapsto \p^G_b = \p_b|_{V^G}$. 
\end{lem}
\proof The action of $\Aut(\G)$ on $E = \oplus_{i\in I}\R\a_i$ respects the bilinear form $(\,.\,|\, .\,)_\G$, and this clearly implies that $g(\p_{i}(e_\a)) = \p_{g(i)}(e_{g(\a)})$ in view of the formulas of definition \ref{def psii}. The result follows by linearity and induction on $\ell(b)$.\qed

\begin{lem}\label{stabilisation fi}
 Assume that $\m f_{i,\a} = \m f_{g(i),g(\a)}$ for every $(i,\a,g) \in I\times \Phi^+ \times G$. Then for every $i,\, j\in I$, $v \in V$ and $g \in G$, we get $f_{g(i)}(g(v)) = f_i(v)$ and $f_{g(i)}\p_{g(j)}(g(v)) = f_i\p_j(v)$. In particular, the linear forms $f_i$ and $f_{g(i)}$ (resp. $f_i\p_j$ and $f_{g(i)}\p_{g(j)}$) coincide on $V^G$.
\end{lem}
\proof The assumption means that $f_{g(i)}(g(e_\a)) = f_i(e_\a)$ for every $(i,\a,g) \in I\times \Phi^+ \times G$, whence the first point by linearity. The second point follows from the first one and the previous lemma. \qed

\begin{pro}\label{stabilisation}Assume that $\m f_{i,\a} = \m f_{g(i),g(\a)}$ for every $(i,\a,g) \in I\times \Phi^+ \times G$. Then for every $(b,v,g) \in B^+\times V \times \Aut(\G)$, we get $g\left(\psi_b(v)\right) = \psi_{g(b)}(g(v))$. In particular, for every $b \in (B^+)^G$, $\psi_b$ stabilizes $V^G$ and hence $\psi$ induces a linear representation $$\psi^G : (B^+)^G \to \mathscr{L}(V^G), \ \ b \mapsto \psi^G_b = \psi_b|_{V^G}.$$

Moreover if the images of $\psi$ are invertible, then so are the images of $\psi^G$. 
\end{pro}
\proof By definition~\ref{def psii}, we get
$\begin{cases}
g(\psi_{i}(e_\a)) = g(\p_{i}(e_\a)) + \m f_{i,\a}e_{\a_{g(i)}}, \text{  and }\\
\psi_{g(i)}(e_{g(\a)}) = \p_{g(i)}(e_{g(\a)}) + \m f_{g(i),g(\a)}e_{\a_{g(i)}}.
\end{cases}$ 

Whence $g(\psi_{i}(e_\a)) = \psi_{g(i)}(e_{g(\a)})$ by assumption and lemma \ref{stabilisation phi}, and the first point by linearity and induction on $\ell(b)$. Moreover if the images of $\psi$ are invertible, that is, if $\m f_{i,\a_i} = f_i(e_{\a_i}) \in \m R^\times$ for every $i \in I$, then the formulas of lemma \ref{condition d'inversibilité} show that we also get $g(\psi_{i}^{-1}(e_\a)) = \psi_{g(i)}^{-1}(e_{g(\a)})$, and hence, similarly to what as just been done, $\psi^{-1}_b$ stabilizes $V^G$ for every $b \in (B^+)^G$. This gives the result. \qed

\begin{defi}
Under the assumption of the previous proposition, we call \emph{twisted LK-representation} the linear representation $\psi^G : (B^+)^G \to \EL(V^G)$ of the Artin-Tits monoid $(B^+)^G = B^+_{\G'}$, and, when appropriate, the induced linear representation $\psi^G_{\gr}: B_{\G'} \to \GL(V^G)$ of the Artin-Tits group $B_{\G'}$.
\end{defi}

The assumption $\m f_{i,\a} = \m f_{g(i),g(\a)}$ for every $(i,\a,g) \in I\times \Phi^+ \times G$ is not always satisfied : for example if $i$ and $g(i)$ are not in the same connected component of $\G$, then $\m f_{i,\a_i}$ and $\m f_{g(i),\a_{g(i)}}$ can be chosen to be distinct (see remark \ref{rels irreductibilite} above). I do not know if this assumption is always satisfied when $\G$ is connected, but we have the following partial result : 

\begin{pro}\label{cas favorables}Let $(f_i)_{i\in I}$ be an LK-family, and assume that we are in one of the following cases : 
\begin{enumerate}
 \item $\G$ is spherical and irreducible (\ie of type $ADE$), or
 \item $\G$ is affine (\ie of type $\tilde A\tilde D\tilde E$), or
 \item $(f_i)_{i\in I}$ is the family constructed in definition \ref{defpolParis} and does not depend on the choice of the $j_\a$'s (for example if $\G$ has no triangle).
 \end{enumerate}

Then $\m f_{i,\a} = \m f_{g(i),g(\a)} \ \, \text{  for every  } \ \, (i,\a,g) \in I\times \Phi^+\times \Aut(\G).$
\end{pro}
\proof The result for the three situations (note that the first one is a consequence of the third one) are easy to see by induction on $\dep(\a)$, using the inductive construction of $(f_i)_{i\in I}$, and the independence results at the inductive steps, of subsection \ref{cas spherique}, \ref{The affine case}, or \ref{The construction of Paris} respectively, and the fact that the action of $\Aut(\G)$ on $E = \oplus_{i\in I}\R\a_i$ respects the bilinear form $(\,.\,|\, .\,)_\G$. \qed

\medskip

We denote by $\Phi^+/G$ the set of orbits of $\Phi^+$ under $G$ and, for every $\Theta \in \Phi^+/G$, we set $e_\Theta = \sum_{\a\in \Theta}e_\a$. The family $(e_\Theta)_{\Theta \in \Phi^+/G}$ is a basis of $V^G$. 

\subsection{Twisted faithfulness criterion}\label{Faithfulness of the induced representation}\mbox{}\medskip

The aim of this subsection is to prove that the faithfulness criterion of subsection \ref{proof of Hee} also works for a twisted LK-representation $\psi^G$. 

\begin{lem}\label{reduction2}Let $\rho : (B^+)^G \to M$ be a monoid homomorphism where $M$ is left cancellative. If $\rho$ satisfies $\rho(b) = \rho(b') \Rightarrow I(b)= I(b')$ for all $b,\,b' \in (B^+)^G$, then it is injective.
\end{lem}
\proof Let $b$, $b' \in (B^+)^G$ be such that $\rho(b) = \rho(b')$. We prove by induction on $\ell(b)$ that $b = b'$. If $\ell(b) = 0$, \ie if $b = 1$, then $I(b)= I(b') = \emptyset$, hence $b' = 1$ and we are done. 
 
If $\ell(b) > 0$, fix $i \in I(b) = I(b')$. Since the action of $G$ on $B^+$ respects the divisibility and since $b$ is fixed by $G$, the orbit $J$ of $i$ under $G$ is included in $I(b) = I(b')$, but then $J$ is spherical and there exist $b_1$, $b'_1\in B^+$ such that $b = \D_Jb_1$ and $b' = \D_Jb'_1$. Since $J$ is an orbit of $I$ under $G$, the element $\D_J$ is fixed by $G$ and hence so are $b_1$ and $b'_1$, so we get $\rho(\D_J)\rho(b_1) = \rho(\D_J)\rho(b'_1)$ in $M$, whence $\rho(b_1) = \rho(b'_1)$ by cancellation, therefore $b_1 = b'_1$ by induction and finally $b = b'$. \qed

\medskip

Let us assume that the condition of proposition \ref{stabilisation} is satisfied, so that the twisted LK-representation $\psi^G : (B^+)^G \to \EL(V^G)$ is defined. Then : 

\begin{thm}\label{main thm}Assume that the following two conditions are satisfied : 
\begin{enumerate}
\item $\Image(\psi^G)$ is a left cancellative submonoid of $\EL(V^G)$, 
\item there exists a totally ordered commutative ring $\m R_0$ and a ring homomorphism $\m R \to \m R_0$, $\m x \mapsto \overline{\m x}$, such that $\overline{\m a}$, $\overline{\m b}$, $\overline{\m c}$, $\overline{\m d}$ are positive and $\overline{\m f_{i,\a}} = 0$ for every $(i,\a) \in I\times \Phi^+$. 
\end{enumerate}

Then the twisted LK-representation $\psi^G$ is faithful.
\end{thm}
\proof Note first that, with notations~\ref{notations bin}, if $R \in \Bin(\Omega)$ and if $(\Psi_\lambda)_{\lambda \in \Lambda}$ is a family of subsets of $\Omega$, we get $R(\bigcup_{\lambda \in \Lambda}\Psi_\lambda) = \bigcup_{\lambda \in \Lambda}R(\Psi_\lambda)$.

In order to prove the theorem, it suffices to show that $\psi^G$ satisfies the assumption of lemma \ref{reduction2}. So let $b,\, b' \in (B^+)^G$ be such that $\psi^G_b = \psi^G_{b'}$ and let us show that $I(b) = I(b')$. Since $\psi_{b}$ and $\psi_{b'}$ coincide on $V^G$, we get in particular $\psi_b(e_\Theta) = \psi_{b'}(e_\Theta)$ for every $\Theta \in \Phi^+/G$. 

With the notations of definition \ref{relations de Digne}, let us consider the set $R_b(\Theta)$. Since the coefficients of the matrix of $\overline{\psi_{b}}$ in the basis $(e_\a)_{\a \in \Phi^+}$ of $V_0$ are non-negative, the set $R_b(\Theta)$ is precisely the set of those indices $\b \in \Phi^+$ for which the coefficient of $e_\b$ in the decomposition of $\overline{\psi_b}(e_\Theta)$ in the basis $(e_\a)_{\a \in \Phi^+}$ is positive. 

The same occurs for $b'$, and hence $\psi_b(e_\Theta) = \psi_{b'}(e_\Theta)$ implies $R_b(\Theta) = R_{b'}(\Theta)$. Since we have $\Phi^+ = \bigcup_{\Theta \in \Phi^+/G}\Theta$, we thus get $R_b(\Phi^+) = \bigcup_{\Theta \in \Phi^+/G}R_b(\Theta) = \bigcup_{\Theta \in \Phi^+/G}R_{b'}(\Theta) = R_{b'}(\Phi^+)$, and hence $I(b) = I(b')$ by lemmas \ref{Hee3} and \ref{Hee2}. \qed

\begin{rem}
Assume that we are in the typical situation of subsections \ref{Comments} and \ref{Comments fin}, so that condition (ii) of theorem \ref{main thm} is satisfied. Then proposition \ref{stabilisation} shows that $\Image(\psi^G)$ is included in $\EL(V^G) \cap \GL((V')^G)$ and hence condition (i) of theorem \ref{main thm} is also satisfied. Hence that twisted LK-representation $\psi^G$ is faithful and so is the induced twisted LK-representation $\psi^G_{\gr} : B_{\G'}\to \GL((V')^G)$ when $\G'$ is spherical. 
\end{rem}

\subsection{Formulas when $|G| = 2$}\label{formules G = 2}\mbox{}\medskip

Recall that $(B^+)^G = B^+_{\G'}$ is generated by the elements $\D_J$, for $J$ running through the spherical orbits of $I$ under $G$ (see section \ref{Graph automorphisms}). In this subsection, we assume that $\m f_{i,\a} = \m f_{g(i),g(\a)}$ for every $(i,\a,g)\in I\times \Phi^+ \times G$, and we compute the maps $\psi^G_{\D_J}$ for those orbits $J$, at least when $|G| = 2$. 

\begin{nota}Let $J$ be an orbit of $I$ under $G$. Then in view of lemma \ref{stabilisation fi} :
\begin{enumerate}
 \item the linear forms $f_j$, for $j \in J$, coincide on $V^G$, and we set $$f_J = f_j|_{V^G} \in (V^G)^\star, \ \text{  for  } \ j\in J,$$
\item if $J = \{i,j\}$, then the forms $f_i\p_j$ and $f_j\p_i$ coincide on $V^G$, and we set $$f'_J = f_i\p_j|_{V^G} = f_j\p_i|_{V^G} \in (V^G)^\star.$$
\end{enumerate} 
\end{nota}

% \begin{rem}
% In the second case of the previous notations, if $J = \{i,j\}$ with $m_{i,j} = 2$, then we simply get $f'_J = \m d f_J$ (by (ii) of definition \ref{def lkp}).
% \end{rem}

Note that if $J$ is an orbit of $I$ under $G$, then $\Theta_J:=\{\a_i\mid i\in J\}$ is an orbit of $\Phi^+$ under $G$ ; moreover if $J = \{i,j\}$ with $m_{i,j} = 3$, then $\Theta'_J = \{\a_i + \a_j\}$ is also an orbit of $\Phi^+$ under $G$.

\begin{pro}\label{forme de psiG}Let $J$ be an orbit of $I$ under $G$. Then  
\begin{enumerate}
\item if $J = \{i\}$, $\psi^G_{\D_J} = \p^G_{\D_J}+f_J\bt e_{\Theta_J}$,
\item if $J = \{i,j\}$ with $m_{i,j} = 2$, $\psi^G_{\D_J} = \p^G_{\D_J} + \m d f_J\bt e_{\Theta_J}$, 
\item if $J = \{i,j\}$ with $m_{i,j} = 3$, $\psi^G_{\D_J} = \p^G_{\D_J} + (\m{bc}f_J+\m af'_J)\bt e_{\Theta_J} + \m cf'_J\bt e_{\Theta'_J}$.
\end{enumerate}
\end{pro}
\proof If $J = \{i\}$, then $\D_J = \gras s_i$ and (i) is clear. If $J = \{i,j\}$ with $m_{i,j} = 2$, then $\D_J = \gras s_i\gras s_j$ and we get $\psi_i\psi_j = \p_i\p_j + f_i\p_j\bt e_{\a_i} + \m df_j\bt e_{\a_j}$ (see the proof of lemma \ref{condition sur les psii}), whence (ii) since $f_i\p_j = \m df_i$. Finally if $J = \{i,j\}$ with $m_{i,j} = 3$, then $\D_J = \gras s_i\gras s_j\gras s_i$ and we get, following the computations of the proof of lemma \ref{condition sur les psii}, $\psi_i\psi_j\psi_i = \p_i\p_j\p_i + f_i\p_j\p_i\bt e_{\a_i} + (\m{bc}f_i + \m a f_j\p_i)\bt e_{\a_j} + \m c f_j\p_i\bt e_{\a_i+\a_j}$ (using the fact that $f_i\p_j(e_{\a_i}) = f_j\p_i(e_{\a_i}) = 0$), whence (iii) since $f_i\p_j\p_i|_{V^G} = f_j\p_i^2|_{V^G}$ and since we have seen, again in the proof of lemma \ref{condition sur les psii}, that $f_j\p_i^2 = \m{bc}f_j + \m a f_j\p_i$.  \qed

\medskip

Hence when $|G| = 2$, the previous proposition gives a complete description of the possible maps $\psi^G_{\D_J}$ when $J$ runs through the (spherical) orbits of $I$ under $G$. 

We detail below the matrix of $\p^G_{\D_J}$ in the basis $(e_\Theta)_{\Theta \in \Phi^+/G}$ and the values on this basis of the linear forms involved in the expression of $\psi^G_{\D_J}$. Note that since the maps $\p_i$, $i\in J$, stabilize the submodule of $V$ generated by the elements $e_\b$ for $\b$ running through a given $J$-mesh $M$, the map $\p^G_J$ stabilizes the submodule of $V^G$ generated by the elements $e_\Theta$ for $\Theta$ running through the orbits included in $G(M)$. The matrix of $\p^G_J$ in the basis $(e_\Theta)_{\Theta \in \Phi^+/G}$ of $V^G$ is then block diagonal, for the corresponding block decomposition.

\subsubsection{Case $J = \{i\}$}\label{J singleton}\mbox{}\medskip 

\noindent $\bullet$ $f_J(e_\Theta) = f_i(e_\Theta) = \begin{cases} \m f_{i,\a} & \text{if  } \Theta = \{\a\} \\ \m f_{i,\a} + \m f_{i,\a'} = 2\m f_{i,\a} & \text{if  } \Theta = \{\a,\a'\}, \a \neq \a'\end{cases}$, 

%\medskip

\noindent $\bullet$ the blocks of $\p^G_{\D_J} = \p_i|_{V^G}$ are the following ones :
\begin{center}
$\left\{\begin{array}{cl}
\overset{\text{\ \scriptsize $e_{\Theta}$}}{\big(\begin{matrix} 0 \end{matrix}\big)}  & \text{if   }\ \ \ \ \ \ \ \ \ \  \Theta = \Theta_J \\ 
\overset{\text{\ \scriptsize $e_{\Theta}$}}{\big(\begin{matrix} \m d\end{matrix}\big)}  & 
\text{if  
\begin{picture}(40,15)(0,-4)
\put(20,0){\circle*{4}}
\put(25,0){\circle{10}}
\put(24,-3){\scriptsize $i$}
\put(0,-3){\small $\Theta$}
\put(15,-4){\dashbox(16,8){}}
\end{picture}  or \
\begin{picture}(55,0)(0,-4)
\put(20,0){\circle*{4}}
\put(25,0){\circle{10}}
\put(24,-3){\scriptsize $i$}
\put(40,0){\circle*{4}}
\put(45,0){\circle{10}}
\put(44,-3){\scriptsize $i$}
\put(0,-3){\small $\Theta$}
\put(15,-4){\dashbox(36,8){}}
\end{picture} } \\ 
\overset{\text{\scriptsize \ \ $e_{\Theta_1}$\ \ $e_{\Theta_2}$}}{\begin{pmatrix}\m a & \m b \\ \m c & 0 \end{pmatrix}} & \text{if 
\begin{picture}(40,25)(0,4)
\put(20,0){\circle*{4}}
\put(20,20){\circle*{4}}
\put(20,0){\line(0,1){20}}
\put(22,8){\scriptsize $i$}
\put(0,-3){\small  $\Theta_1$}
\put(0,17){\small $\Theta_2$}
\put(15,-4){\dashbox(10,8){}}
\put(15,16){\dashbox(10,8){}}
\end{picture}  or \
\begin{picture}(55,25)(0,4)
\put(20,0){\circle*{4}}
\put(20,20){\circle*{4}}
\put(20,0){\line(0,1){20}}
\put(22,8){\scriptsize $i$}
\put(0,-3){\small  $\Theta_1$}
\put(0,17){\small $\Theta_2$}
\put(15,-4){\dashbox(30,8){}}
\put(15,16){\dashbox(30,8){}}
\put(40,0){\circle*{4}}
\put(40,20){\circle*{4}}
\put(40,0){\line(0,1){20}}
\put(42,8){\scriptsize $i$}
\end{picture}}
\end{array}\right.$
\end{center}

%\medskip
%\vspace*{2mm}

\subsubsection{Case $J = \{i,j\}$ with $m_{i,j} = 2$}\label{J avec mij=2}\mbox{}\medskip 

\noindent $\bullet$ $f_J(e_\Theta) = \m {df}_i(e_\Theta) = \begin{cases} \m {df}_{i,\a} & \text{if  } \Theta = \{\a\} \\ \m d(\m f_{i,\a} + \m f_{i,\a'}) = \m d(\m f_{i,\a}+\m f_{j,\a}) & \text{if  } \Theta = \{\a,\a'\}, \a \neq \a'\end{cases}$, 

%\medskip

\noindent $\bullet$ the blocks of $\p^G_{\D_J} = (\p_i\p_j)|_{V^G}$ are the following ones :
\begin{center}
$\left\{\begin{array}{cl}
  \overset{\text{\ \scriptsize $e_{\Theta}$}}{\big(\begin{matrix} 0 \end{matrix}\big)} & \text{if \  \ \ \ \ \ \ \ \ \ \ \ $\Theta = \Theta_J$} \\
\overset{\text{\ \scriptsize $e_{\Theta}$}}{\big(\begin{matrix} \m{d}^2\end{matrix}\big)} & \text{if \ \
\begin{picture}(40,18)(0,-4)
\put(20,0){\circle*{4}}
\put(25,0){\circle{10}}
\put(24,-3){\scriptsize $i$}
\put(15,0){\circle{10}}
\put(12,-3){\scriptsize $j$}
\put(-3,-3){\small $\Theta$}
\put(8,-4){\dashbox(24,8){}}
\end{picture} or \ 
\begin{picture}(60,10)(0,-4)
\put(20,0){\circle*{4}}
\put(25,0){\circle{10}}
\put(24,-3){\scriptsize $i$}
\put(15,0){\circle{10}}
\put(12,-3){\scriptsize $j$}
\put(45,0){\circle*{4}}
\put(50,0){\circle{10}}
\put(49,-3){\scriptsize $i$}
\put(40,0){\circle{10}}
\put(37,-3){\scriptsize $j$}
\put(-3,-3){\small $\Theta$}
\put(8,-4){\dashbox(50,8){}}
\end{picture}} \\ 
\overset{\text{\scriptsize \ \ $e_{\Theta_1}$\ \ $e_{\Theta_2}$}}{\begin{pmatrix}\m a \m d & \m b \m d \\ \m c \m d & 0 \end{pmatrix}} & \text{if \ \ 
\begin{picture}(55,30)(-20,4)
\put(20,0){\circle*{4}}
\put(20,20){\circle*{4}}
\put(20,0){\line(0,1){20}}
\put(22,8){\scriptsize $i$}
\put(15,0){\circle{10}}
\put(12,-3){\scriptsize $j$}
\put(15,20){\circle{10}}
\put(12,17){\scriptsize $j$}
\put(-6,-3){\small  $\Theta_1$}
\put(-6,17){\small $\Theta_2$}
\put(8,-4){\dashbox(36,8){}}
\put(8,16){\dashbox(36,8){}}
\put(40,0){\circle*{4}}
\put(40,20){\circle*{4}}
\put(40,0){\line(0,1){20}}
\put(42,8){\scriptsize $j$}
\put(35,0){\circle{10}}
\put(32,-3){\scriptsize $i$}
\put(35,20){\circle{10}}
\put(32,17){\scriptsize $i$}
\end{picture}}\\ 
\overset{\text{\scriptsize \ \ \ $e_{\Theta_1}$\ \ \ \ $e_{\Theta_2}$ \ \ $e_{\Theta_3}$}}{\begin{pmatrix}
\m a^2 & 2\m {ab} & \m b^2 \\
\m {ac} & \m {bc} & 0 \\
\m c^2 & 0 & 0
\end{pmatrix}}  & \text{if \ \ 
\begin{picture}(80,40)(-30,14)
\put(20,0){\circle*{4}}
\put(0,20){\circle*{4}}
\put(40,20){\circle*{4}}
\put(20,40){\circle*{4}}
\put(20,0){\line(1,1){20}}
\put(20,0){\line(-1,1){20}}
\put(20,40){\line(1,-1){20}}
\put(20,40){\line(-1,-1){20}}
\put(33,5){\scriptsize $i$}
\put(33,30){\scriptsize $j$}
\put(2,5){\scriptsize $j$}
\put(5,30){\scriptsize $i$}
\put(-23,-3){\small  $\Theta_1$}
\put(-23,17){\small $\Theta_2$}
\put(-23,37){\small $\Theta_3$}
\put(15,-4){\dashbox(10,8){}}
\put(15,36){\dashbox(10,8){}}
\put(-5,16){\dashbox(50,8){}}
\end{picture}} \\ 
\overset{\text{\scriptsize \ \ $e_{\Theta_1}$\ \ $e_{\Theta_2}$ \ \ $e_{\Theta_3}$  \  $e_{\Theta_4}$}}{\begin{pmatrix}
\m a^2 & \m {ab} & \m {ab} & \m b^2 \\
\m {ac} & 0 &\m {bc} & 0\\
\m {ac} & \m {bc} & 0 & 0 \\
\m c^2 & 0 & 0 & 0
\end{pmatrix}} & \text{if \ \ 
\begin{picture}(100,50)(-20,14)
\put(20,0){\circle*{4}}
\put(0,20){\circle*{4}}
\put(40,20){\circle*{4}}
\put(20,40){\circle*{4}}
\put(20,0){\line(1,1){20}}
\put(20,0){\line(-1,1){20}}
\put(20,40){\line(1,-1){20}}
\put(20,40){\line(-1,-1){20}}
\put(33,8){\scriptsize $i$}
\put(33,29){\scriptsize $j$}
\put(3,7){\scriptsize $j$}
\put(5,30){\scriptsize $i$}
\put(-23,-3){\small  $\Theta_1$}
\put(-23,17){\small $\Theta_2$}
\put(-23,37){\small $\Theta_4$}
\put(15,-4){\dashbox(30,8){}}
\put(15,36){\dashbox(30,8){}}
\put(-5,16){\dashbox(30,8){}}
\put(40,0){\circle*{4}}
\put(20,20){\circle*{4}}
\put(60,20){\circle*{4}}
\put(40,40){\circle*{4}}
\put(40,0){\line(1,1){20}}
\put(40,0){\line(-1,1){20}}
\put(40,40){\line(1,-1){20}}
\put(40,40){\line(-1,-1){20}}
\put(53,7){\scriptsize $j$}
\put(52,30){\scriptsize $i$}
\put(35,16){\dashbox(30,8){}}
\put(73,17){\small $\Theta_3$}
\end{picture}}
 \end{array}\right.$
\end{center}

%\medskip
%\vspace*{2mm}

\subsubsection{Case $J = \{i,j\}$ with $m_{i,j} = 3$}\label{J avec mij=3}\mbox{}\medskip

We detail below the values of the linear forms $\m{bc}f_{J}+\m{a}f'_J$ and $\m cf'_{J}$, and the blocks of $\p^G_{\D_J} = (\p_i\p_j\p_i)|_{V^G}$ for the different possible configurations of orbits. 

%\medskip 

\noindent $\bullet$ Orbits $\Theta_J = \{\a_i,\a_j\}$ and $\Theta'_J = \{\a_i+\a_j\}$.
\begin{center}
\begin{tabular}{|c|c|c||c|}\hline Orbits & Values of $\m{bc}f_{J}+\m{a}f'_J$ &  Values of $\m cf'_{J}$ & Block in  $\p^{\overset{\mbox{}}{G}}_{\D_J}$ \\ \hline \hline 
\begin{tabular}{c}$\Theta_J$ \\ $\Theta'_J$ \end{tabular} & \begin{tabular}{c}$\m{bcf}_{i,\a_i}$ \\ 0 \end{tabular} & \begin{tabular}{c} 0  \\ $\m{bcf}_{i,\a_i}$ \end{tabular}  & $\begin{pmatrix}0 & 0\\ 0 & 0 \end{pmatrix}$\\ \hline
\end{tabular}
\end{center}

%\medskip

\noindent$\bullet$ Configuration \ 
\begin{picture}(30,15)(0,-2)
\put(20,0){\circle*{4}}\put(18,6){$\a$}
\put(25,0){\circle{10}}
\put(24,-3){\scriptsize $i$}
\put(15,0){\circle{10}}
\put(12,-3){\scriptsize $j$}
\put(-3,-3){\small $\Theta$}
\put(8,-4){\dashbox(24,8){}}
\end{picture} \ \  or\ \ \begin{picture}(60,15)(0,-2)
\put(20,0){\circle*{4}}\put(18,6){$\a$}
\put(25,0){\circle{10}}
\put(24,-3){\scriptsize $i$}
\put(15,0){\circle{10}}
\put(12,-3){\scriptsize $j$}
\put(45,0){\circle*{4}}\put(43,6){$\a'$}
\put(50,0){\circle{10}}
\put(49,-3){\scriptsize $i$}
\put(40,0){\circle{10}}
\put(37,-3){\scriptsize $j$}
\put(-3,-3){\small $\Theta$}
\put(8,-4){\dashbox(50,8){}}
\end{picture}

\begin{center}
\begin{tabular}{|c|c|c|}\hline Value of $\m{bc}f_{J}+\m{a}f'_J$  &  Value of $\m cf'_{J}$ & Block in  $\p^{\overset{\mbox{}}{G}}_{\D_J}$ \\ \hline \hline 
 $|\Theta|\m d^2 \m f_{i,\a}$ &  $|\Theta|\m{cd} \m f_{i,\a}$ & $\overset{\text{\ \scriptsize $e_{\Theta}$}}{\big(\begin{matrix}\m d^3\end{matrix}\big)} $ \\ \hline
\end{tabular}
\end{center}

\noindent$\bullet$ Configuration \ \ \begin{picture}(70,35)(0,16)
%\put(20,0){\circle*{4}}
\put(20,0){\circle*{4}}
\put(20,20){\circle*{4}}
\put(20,40){\circle*{4}}
\put(20,0){\line(0,1){40}}
\put(22,8){\scriptsize $i$}
\put(22,28){\scriptsize $j$}
\put(15,0){\circle{10}}
\put(12,-2){\scriptsize $j$}
\put(15,40){\circle{10}}
\put(12,37){\scriptsize $i$}
\put(-5,-3){\small  $\Theta_1$}
\put(-5,17){\small $\Theta_2$}
\put(-5,37){\small $\Theta_3$}
\put(48,-3){$\g$}
\put(48,17){$\b$}
\put(48,37){$\a$}
\put(8,-4){\dashbox(36,8){}}
\put(8,36){\dashbox(36,8){}}
\put(15,16){\dashbox(30,8){}}
\put(40,0){\circle*{4}}
\put(40,20){\circle*{4}}
\put(40,40){\circle*{4}}
\put(40,0){\line(0,1){40}}
\put(42,8){\scriptsize $j$}
\put(35,0){\circle{10}}
\put(32,-3){\scriptsize $i$}
\put(42,28){\scriptsize $i$}
\put(35,40){\circle{10}}
\put(32,38){\scriptsize $j$}
\end{picture}

%\vspace{8mm}

\begin{center}
\begin{tabular}{|c|c|c||c|}\hline Orbits & Values of $\m{bc}f_{J}+\m{a}f'_J$ &  Values of $\m cf'_{J}$ & Block in  $\p^{\overset{\mbox{}}{G}}_{\D_J}$ \\ \hline \hline 
\begin{tabular}{c}$\Theta_1$ \\ $\Theta_2$ \\ $\Theta_3$ \end{tabular} & 
$\begin{array}{l}
  \m{bcf}_{i,\g}+\m{d}(\m d +\m a)\m f_{j,\g}\\
 \m b\left(\m{cf}_{j,\b}+\m{af}_{i,\g}+\m{df}_{j,\g}\right) \\
 \m b\left(\m{df}_{j,\b} + \m{bf}_{i,\g}\right) 
\end{array}$
& $\begin{array}{l}
  2\m{cdf}_{j,\g}\\
 2\m{bcf}_{i,\g} \\
 2\m{bcf}_{j,\b}
\end{array}$ & $\overset{\text{\scriptsize $\overset{\mbox{}}{e_{\Theta_1}}$\ \ \ $e_{\Theta_2}$ \ \ \  $e_{\Theta_3}$}}{\begin{pmatrix}
\m{ad}^2 & \m{abd} & \m b^2 \m d \\
\m{acd} & \m{bcd} & 0 \\
\m c^2\m d & 0 & 0
\end{pmatrix}}$\\ \hline 
\end{tabular}
\end{center}

\noindent$\bullet$ Configuration \ \ \begin{picture}(70,45)(-20,22)
\put(20,0){\circle*{4}}
\put(0,20){\circle*{4}}
\put(40,20){\circle*{4}}
\put(0,40){\circle*{4}}
\put(40,40){\circle*{4}}
\put(20,60){\circle*{4}}
\put(20,0){\line(1,1){20}}
\put(20,0){\line(-1,1){20}}
\put(0,20){\line(0,1){20}}
\put(40,20){\line(0,1){20}}
\put(20,60){\line(1,-1){20}}
\put(20,60){\line(-1,-1){20}}
\put(33,5){\scriptsize $i$}
\put(33,50){\scriptsize $i$}
\put(43,28){\scriptsize $j$}
\put(-5,28){\scriptsize $i$}
\put(4,5){\scriptsize $j$}
\put(5,50){\scriptsize $j$}
\put(-23,-3){\small  $\Theta_1$}
\put(-23,17){\small $\Theta_2$}
\put(-23,37){\small $\Theta_3$}
\put(-23,57){\small $\Theta_4$}
\put(28,-5){$\delta$}
\put(48,17){$\g$}
\put(48,37){$\b$}
\put(28,57){$\a$}
\put(15,-4){\dashbox(10,8){}}
\put(15,56){\dashbox(10,8){}}
\put(-5,36){\dashbox(50,8){}}
\put(-5,16){\dashbox(50,8){}}
\end{picture} 

\vspace{8mm}

\begin{center}
\begin{tabular}{|c|c|c|}\hline Orbits & Values of $\overset{\mbox{}}{\m{bc}f_{J}+\m{a}f'_J}$ &  Values of $\m cf'_{J}$ \\ \hline \hline 
\begin{tabular}{c}$\overset{\mbox{}}{\Theta_1}$ \\ $\Theta_2$ \\ $\Theta_3$ \\ $\Theta_4$ \end{tabular} &  
$\begin{array}{l}
\overset{\mbox{}}{\m{acf}_{j,\g}}+(\m a^2 + \m{bc})\m f_{i,\delta} \\
\m b\left(\m c(\m f_{i,\g}+\m f_{j,\g})+2\m{af}_{i,\delta}\right) \\
\m b\left(\m{cf}_{j,\b}+\m{af}_{i,\g}+\m{bf}_{i,\delta}\right) \\
\m b^2\m f_{i,\g} \\
\end{array}$ & 
$\begin{array}{ll}
\overset{\mbox{}}{\m c^2\m f_{j,\g}}+\m{acf}_{j,\delta} \\
2\m{bcf}_{i,\delta} \\
2\m{bcf}_{i,\g} \\
\m{bcf}_{j,\b} \\
\end{array}$ \\ \hline \hline 
\multicolumn{3}{|c|}{Block in  $\p^{\overset{\mbox{}}{G}}_{\D_J}$} \\ \hline \hline
\multicolumn{3}{|c|}{
$\overset{\text{\scriptsize \ \ \ \ \  \ \  $\overset{\mbox{}}{e_{\Theta_1}}$ \ \ \ \ \ \ \ \ \ $e_{\Theta_2}$  \ \ \ \ \   $e_{\Theta_3}$  \ \ \   $e_{\Theta_4}$}}{\begin{pmatrix}
\m a(\m a^2 + \m{bc}) & 2\m a^2\m b & 2\m{ab}^2 & \m b^3 \\
\m a^2\m c & 2\m{abc} & \m b^2\m c & 0\\
\m{ac}^2 & \m{bc}^2 & 0 & 0 \\
\m{c}^3 & 0 & 0 & 0
\end{pmatrix}}$} \\ \hline 
\end{tabular}
\end{center}

\noindent$\bullet$ Configuration \ \ \begin{picture}(110,55)(-20,20)
\put(20,0){\circle*{4}}
\put(0,20){\circle*{4}}
\put(40,20){\circle*{4}}
\put(0,40){\circle*{4}}
\put(40,40){\circle*{4}}
\put(20,60){\circle*{4}}
\put(20,0){\line(1,1){20}}
\put(20,0){\line(-1,1){20}}
\put(0,20){\line(0,1){20}}
\put(40,20){\line(0,1){20}}
\put(20,60){\line(1,-1){20}}
\put(20,60){\line(-1,-1){20}}
\put(33,8){\scriptsize $j$}
\put(33,49){\scriptsize $j$}
\put(33,27){\scriptsize $i$}
\put(23,27){\scriptsize $i$}
\put(3,7){\scriptsize $i$}
\put(3,50){\scriptsize $i$}
\put(-4,28){\scriptsize $j$}
\put(-23,-3){\small  $\Theta_1$}
\put(-23,17){\small $\Theta_2$}
\put(-23,37){\small $\Theta_4$}
\put(-23,57){\small $\Theta_6$}
\put(8,66){$\a$}
\put(42,66){$\a'$}
\put(-10,47){$\b$}
\put(64,47){$\b'$}
\put(-10,10){$\g$}
\put(64,10){$\g'$}
\put(8,-10){$\delta$}
\put(46,-10){$\delta'$}
\put(15,-4){\dashbox(30,8){}}
\put(15,56){\dashbox(30,8){}}
\put(-5,36){\dashbox(30,8){}}
\put(35,36){\dashbox(30,8){}}
\put(-5,16){\dashbox(30,8){}}
\put(40,0){\circle*{4}}
\put(20,20){\circle*{4}}
\put(60,20){\circle*{4}}
\put(20,40){\circle*{4}}
\put(60,40){\circle*{4}}
\put(40,60){\circle*{4}}
\put(40,0){\line(1,1){20}}
\put(40,0){\line(-1,1){20}}
\put(20,20){\line(0,1){20}}
\put(60,20){\line(0,1){20}}
\put(40,60){\line(1,-1){20}}
\put(40,60){\line(-1,-1){20}}
\put(53,7){\scriptsize $i$}
\put(53,50){\scriptsize $i$}
\put(62,28){\scriptsize $j$}
\put(35,16){\dashbox(30,8){}}
\put(73,17){\small $\Theta_3$}
\put(73,37){\small $\Theta_5$}
\end{picture}

\vspace{8mm}

\begin{center}
\begin{tabular}{|c|c|c|}\hline Orbits & Values of $\overset{\mbox{}}{\m{bc}f_{J}+\m{a}f'_J}$ &  Values of $\m cf'_{J}$  \\ \hline \hline 
\begin{tabular}{c}$\overset{\mbox{}}{\Theta_1}$ \\ $\Theta_2$ \\ $\Theta_3$ \\ $\Theta_4$ \\ $\Theta_5$ \\ $\Theta_6$ \end{tabular} & 
$\begin{array}{l}
\m{bc}\left(\m f_{i,\delta}+\m f_{j,\delta}\right) + \overset{\mbox{}}{2\m{a}}\left(\m{af}_{j,\delta}+\m{cf}_{j,\g}\right) \\
\m b\left(\m c(\m f_{i,\g} + \m f_{j,\g}) + 2\m{af}_{i,\delta'}\right) \\ 
\m b\left(\m c(\m f_{i,\g'} + \m f_{j,\g'}) + 2\m{af}_{i,\delta}\right) \\ 
\m b\left(\m{cf}_{j,\b}+\m{af}_{i,\g}+\m{bf}_{j,\delta}\right)\\
\m b\left(\m{cf}_{j,\b'}+\m{af}_{i,\g'}+\m{bf}_{j,\delta'}\right)\\
\m b^2(\m f_{i,\g}+\m f_{i,\g'})
\end{array}$ & 
$\begin{array}{l}
\overset{\mbox{}}{2\m c(\m{af}_{j,\delta}}+\m{cf}_{j,\g}) \\ 
2\m{bcf}_{i,\delta'} \\ 
2\m{bcf}_{i,\delta} \\ 
2\m{bcf}_{i,\g}\\
2\m{bcf}_{i,\g'}\\
2\m{bcf}_{j,\b}
\end{array}$ \\ \hline \hline 
\multicolumn{3}{|c|}{Block in  $\p^{\overset{\mbox{}}{G}}_{\D_J}$} \\ \hline \hline
\multicolumn{3}{|c|}{
$\overset{\text{\scriptsize \ \ \ \ \ \ \ $\overset{\mbox{}}{e_{\Theta_1}}$ \ \ \ \ \ \ \ \ $e_{\Theta_2}$ \ \ \ \ $e_{\Theta_3}$ \ \ \ $e_{\Theta_4}$  \ \ \  $e_{\Theta_5}$ \ \  $e_{\Theta_6}$}}{\begin{pmatrix}
\m a(\m a^2 + \m{bc}) & \m a^2\m b & \m a^2\m b & \m{ab}^2 & \m{ab}^2 & \m b^3 \\
\m a^2\m c & \m{abc} & \m{abc} & 0 & \m b^2\m c & 0\\
\m a^2\m c & \m{abc} & \m{abc} & \m b^2\m c & 0 & 0\\
\m{ac}^2 & 0 & \m{bc}^2 & 0 & 0 & 0 \\
\m{ac}^2 & \m{bc}^2 & 0 & 0 & 0 & 0\\
\m{c}^3 & 0 & 0 & 0 & 0 & 0
\end{pmatrix}} $ }\\ \hline 
\end{tabular}
\end{center}

%\medskip

\subsection{Twisted LK-representations of type $B$}\label{Twisted LK-representations of type B}\mbox{}\medskip

Fix $n \in \N_{\geqslant 3}$. Then the Artin-Tits group $B$ of Coxeter type $B_n$ appears as the subgroup of fixed elements of three Artin-Tits groups of type $\G_1 = A_{2n-1}$, $\G_2 = A_{2n}$ and $\G_3 = D_{n+1}$ respectively, under a group of graph automorphisms $G_1$, $G_2$ and $G_3$ respectively, where $G_k = \Aut(\G_k)$ for $1\leqslant k \leqslant 3$, except for $n = k = 3$ where $G_3$ is a subgroup of order two of $\Aut(D_4)$ (see \cite{C2,Ca}). We denote by $I_k$ the vertex set of $\G_k$ for $1\leqslant k \leqslant 3$.

%\medskip

Now fix a commutative ring $\m R$ and $\m{(b,c,d)} \in (\m R^\times)^3$, and consider three LK-representations $\psi_1 : B_{\G_1} \to \GL(V_1)$, $\psi_2 : B_{\G_2} \to \GL(V_2)$ and $\psi_3 : B_{\G_3} \to \GL(V_3)$. Recall that $\psi_k$, $1\leqslant k \leqslant 3$, is determined by the common value $\m f_k \in \m R^\times$ of the $\m f_{i,\a_i} = f_i(e_{\a_i})$ for $i \in I_k$ (see subsection \ref{cas spherique}). In view of propositions \ref{stabilisation} and \ref{cas favorables} above, we get three twisted LK-representations $\psi_1^{G_1} : B \to \GL(V_1^{G_1})$,  $\psi_2^{G_2} : B \to \GL(V_2^{G_2})$ and $\psi_3^{G_3} : B \to \GL(V_3^{G_3})$ of the Artin-Tits group $B$. Note that $\psi_1^{G_1}$ is essentially the representation of $B$ considered in~\cite{Di}. 

\medskip

The representation $\psi_2^{G_2}$ is trivially non-equivalent to the two others since it is of degree $|\Phi^+_{\G_2}/G_2| = n(n+1)$ whereas the two others are of degree $|\Phi^+_{\G_1}/G_1| = |\Phi^+_{\G_3}/G_3| = n^2$. The aim of this section is to show that $\psi_1^{G_1}$ and $\psi_3^{G_3}$ are non-equivalent, when $\m R = \Z[x^{\pm 1}, y^{\pm 1}]$ and $(\m{b,c,d}) = (y^p,y^q,y^r)$ with $p,\, q,\, r \in \Z$ such that $2r \neq p+q$ (as in subsection \ref{Comments}), and at least for all $n\geqslant 3$ but two. 

\begin{nota}Following \cite{Bo}, we label by $1$, $2$, $\ldots, n$, the vertices of the Coxeter graph $B_n$, in such a way that the vertex $n$ is the terminal vertex of the edge labeled $4$, and we denote by $\D_1, \,\ldots,\, \D_n$ the corresponding standard generators of $B$. 

Note that we will keep the same symbols for the standard generators of $B$ when considering this group as an abstract Artin-Tits group, or as the subgroup of fixed elements $(B_{\G_k})^{G_k}$ of $B_{\G_k}$ for $1\leqslant k \leqslant 3$. The meaning of $\D_i$, $1\leqslant i \leqslant n$, in terms of a product of the standard generators of $B_{\G_k}$, $1\leqslant k \leqslant 3$, is given in the following table (where the vertices of $\G_k$ are labeled as in \cite{Bo}) : 
\begin{center}
\begin{tabular}{|c||c|c|} \hline & $1\leqslant i < n$ & $i = n$ \\ \hline \hline 
 $k = 1$ & $\gras s_i \gras s_{2n-i}$ & $\gras s_n$ \\ \hline 
$k = 2$ & $\gras s_i \gras s_{2n+1-i}$ & $\gras s_n\gras s_{n+1}\gras s_n$ \\ \hline
$k = 3$ & $\gras s_i$ & $\gras s_n\gras s_{n+1}$ \\ \hline
\end{tabular}
\end{center}
\end{nota}

\begin{lem}The determinant of the map $(\psi^{G_k}_k)_{\D_i}$, for $1 \leqslant i \leqslant n$ and $1\leqslant k \leqslant 3$, is given in the following table :
\begin{center}
\begin{tabular}{|c||c|c|} \hline & $1\leqslant i < n$ & $i = n$ \\ \hline \hline 
 $k = 1$ & $-(\m{bc})^{\overset{\mbox{}}{2n-1}}\m d^{2n(n-2)}\m f_1$ & $(-1)^{n-1}(\m{bc})^{n-1}\m d^{(n-1)^{\overset{\mbox{}}{2}}}\m f_1$ \\ \hline 
$k = 2$ & $(\m{bc})^{\overset{\mbox{}}{2n}}\m d^{2(n^{\overset{\mbox{}}{2}}-n-1)}\m f_2$ & $(-1)^{n-1}(\m{bc})^{3n-1}\m d^{3n(n-1)}\m f_2^2$ \\ \hline
$k = 3$ & $-(\m{bc})^{2n-3}\m d^{(n-1)^{\overset{\mbox{}}{2}}}\m f_3$ & $(-1)^{n-1}(\m{bc})^{3(n-1)}\m d^{2(n-1)(n-2)}\m f_3$ \\ \hline
\end{tabular}
\end{center}
\end{lem}
\proof In view of the formulas of subsection \ref{formules G = 2}, the determinant of $(\psi^{G_k}_k)_{\D_i}$ is of the form $\det(M)\m f_k$ (resp. $\det(M)(\m{bcf}_2)^2$) if $(i,k) \neq (n,2)$ (resp. $(i,k) = (n,2)$), where $M$ is a block diagonal matrix, with blocks of determinant $\m d$, $-\m{bc}$, $\m d^2$, $-\m{bcd}^2$, $-(\m{bc})^3$ or $(\m{bc})^4$ (resp. $\m d^3$, $-(\m{bcd})^3$, $(\m{bc})^6$ or $-(\m{bc})^9$) depending on the configuration of the corresponding orbit in $\Phi^+_{\G_k}$. The result then follows from a direct computation of the number of occurrences of each configuration in $\Phi^+_{\G_k}$ for $1 \leqslant k \leqslant 3$. \qed

\begin{pro}Assume that $\m R = \Z[x^{\pm 1}, y^{\pm 1}]$ and that $(\m{b,c,d}) = (y^p,y^q,y^r)$, with $p,\, q,\, r \in \Z$ such that $2r \neq p+q$. Then the twisted LK-representations $\psi_1^{G_1}$ and $\psi_3^{G_3}$ are non-equivalent, except possibly for two values of $n$ when $r<0<p+q$ or $p+q < 0 < r$.
\end{pro}
\proof If $\psi_1^{G_1}$ and $\psi_3^{G_3}$ were equivalent, then the determinants of $(\psi^{G_1}_1)_{\D_i}$ and $(\psi^{G_3}_3)_{\D_i}$ for $1 \leqslant i < n$ (resp. of $(\psi^{G_1}_1)_{\D_n}$ and $(\psi^{G_3}_3)_{\D_n}$) should be equal. This would imply, in view of the previous lemma, 
\begin{center}
$\begin{cases}
 (\m{bc})^2 \m d^{n^2-2n-1}\m f_1 = \m f_3, \text{  and} \\ \m f_1 = (\m{bc})^{2(n-1)} \m d^{(n-1)(n-3)}\m f_3
 \end{cases}$, whence $(\m{bc})^{2n} \m d^{2(n^2-3n+1)} = 1$,
\end{center}
and by choice of $(\m{b,c,d})$, this is equivalent to $2n(p+q) + 2(n^2-3n+1)r = 0$. 

It is clear that there is at most two values of $n$ satisfying this equality, and that in such a case, $r$ and $p+q$ cannot be zero or of the same sign (note that $n^2-3n+1 = (n-1)(n-2) - 1 \geqslant 1$ since $n \geqslant 3$). \qed

\subsection{Final remark on $\Phi^+/G$}\label{Final remark}\mbox{}\medskip

Recall that we denote by $\G'$ the type of $W^G$ and $(B^+)^G$.

\medskip

When $\G$ is spherical, it is possible to index the basis $(e_\Theta)_{\Theta \in \Phi^+/G}$ of $V^G$ with the set of positive roots of a {\it finite crystallographic} root system (\ie a root system in the sense of \cite[Ch.~VI]{Bo}) of Weyl group $W_{\G'}$, via the bijection $\Theta \mapsto \a_\Theta$, where $\a_\Theta = \dfrac{1}{\text{Card}(\Theta)}\sum_{\a \in \Theta} \a$ (see \cite[Ch. 13]{Car} for justifications).

For example if $\G = A_{2n-1}$, $A_{2n}$ or $D_{n+1}$ and $G = \Aut(\G)$ (or a subgroup of order $2$ of $\Aut(\G)$ for $D_4$), we get a finite crystallographic root system of Dynkin type $C_n$, $BC_n$ or $B_n$ respectively. 

This change of index set increases the resemblance between the twisted and non-twisted cases, and has been used by Digne in \cite{Di} for his proof of faithfulness of $\psi^G$, in the particular cases $\G = A_{2n-1}$, $E_6$ or $D_4$ and $G = \Aut(\G)$.

\medskip

But this change of index set is not possible in general. Indeed, the map $\Theta \mapsto \a_\Theta$ is not necessarily injective if $\G$ is not spherical : for example when $|G| = 2$, then for the following configurations of orbits 
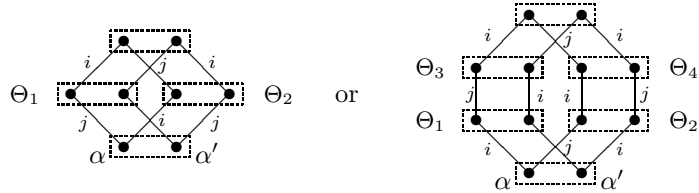
\begin{center}
\begin{picture}(150,50)(-20,-10)
\put(20,0){\circle*{4}}
\put(0,20){\circle*{4}}
\put(40,20){\circle*{4}}
\put(20,40){\circle*{4}}
\put(20,0){\line(1,1){20}}
\put(20,0){\line(-1,1){20}}
\put(20,40){\line(1,-1){20}}
\put(20,40){\line(-1,-1){20}}
\put(33,8){\scriptsize $i$}
\put(33,29){\scriptsize $j$}
\put(3,7){\scriptsize $j$}
\put(5,30){\scriptsize $i$}
%\put(-23,-3){\small  $\Theta$}
\put(-23,17){\small $\Theta_1$}
%\put(-23,37){\small $\Theta'''$}
\put(15,-4){\dashbox(30,8){}}
\put(15,36){\dashbox(30,8){}}
\put(-5,16){\dashbox(30,8){}}
\put(40,0){\circle*{4}}
\put(20,20){\circle*{4}}
\put(60,20){\circle*{4}}
\put(40,40){\circle*{4}}
\put(40,0){\line(1,1){20}}
\put(40,0){\line(-1,1){20}}
\put(40,40){\line(1,-1){20}}
\put(40,40){\line(-1,-1){20}}
\put(53,7){\scriptsize $j$}
\put(52,30){\scriptsize $i$}
\put(35,16){\dashbox(30,8){}}
\put(73,17){\small $\Theta_2$}\put(100,17){\text{or}}
\put(7,-6){$\a$}
\put(47,-6){$\a'$}
\end{picture}
 \begin{picture}(110,60)(-20,0)
\put(20,0){\circle*{4}}
\put(0,20){\circle*{4}}
\put(40,20){\circle*{4}}
\put(0,40){\circle*{4}}
\put(40,40){\circle*{4}}
\put(20,60){\circle*{4}}
\put(20,0){\line(1,1){20}}
\put(20,0){\line(-1,1){20}}
\put(0,20){\line(0,1){20}}
\put(40,20){\line(0,1){20}}
\put(20,60){\line(1,-1){20}}
\put(20,60){\line(-1,-1){20}}
\put(33,8){\scriptsize $j$}
\put(33,49){\scriptsize $j$}
\put(33,27){\scriptsize $i$}
\put(23,27){\scriptsize $i$}
\put(3,7){\scriptsize $i$}
\put(3,50){\scriptsize $i$}
\put(-4,28){\scriptsize $j$}
%\put(-23,-3){\small  $\Theta$}
\put(-23,17){\small $\Theta_1$}
\put(-23,37){\small $\Theta_3$}
%\put(-23,57){\small $\Theta^{(5)}$}
\put(7,-6){$\a$}
\put(47,-6){$\a'$}
% \put(-10,47){$\b_1$}
% \put(64,47){$\b_2$}
% \put(-10,10){$\g_1$}
% \put(64,10){$\g_2$}
% \put(8,-10){$\delta_1$}
% \put(46,-10){$\delta_2$}
\put(15,-4){\dashbox(30,8){}}
\put(15,56){\dashbox(30,8){}}
\put(-5,36){\dashbox(30,8){}}
\put(35,36){\dashbox(30,8){}}
\put(-5,16){\dashbox(30,8){}}
\put(40,0){\circle*{4}}
\put(20,20){\circle*{4}}
\put(60,20){\circle*{4}}
\put(20,40){\circle*{4}}
\put(60,40){\circle*{4}}
\put(40,60){\circle*{4}}
\put(40,0){\line(1,1){20}}
\put(40,0){\line(-1,1){20}}
\put(20,20){\line(0,1){20}}
\put(60,20){\line(0,1){20}}
\put(40,60){\line(1,-1){20}}
\put(40,60){\line(-1,-1){20}}
\put(53,7){\scriptsize $i$}
\put(53,50){\scriptsize $i$}
\put(62,28){\scriptsize $j$}
\put(35,16){\dashbox(30,8){}}
\put(73,17){\small $\Theta_2$}
\put(73,37){\small $\Theta_4$}
\end{picture}
\end{center}
we get $\a_{\Theta_1} = \a_{\Theta_2}$ and $\a_{\Theta_3} = \a_{\Theta_4}$ as soon as $(\a_i|\a) = (\a_j|\a)$.

Note that the first of those counterexamples occurs for example in a root system of type $\tilde A_{2n-1}$ ($n \geqslant 2$) with $G$ generated by the ``half turn'', and the second (which does not occur in the affine cases in view of remark \ref{hexagones affines} above) occurs for example in the root system associated with the Coxeter graph \begin{picture}(20,7)(-5,2)\put(0,0){\circle*{2}}\put(10,0){\circle*{2}}\put(0,10){\circle*{2}}\put(10,10){\circle*{2}}\put(0,0){\line(0,1){10}}\put(0,0){\line(1,0){10}}\put(10,10){\line(0,-1){10}}\put(10,10){\line(-1,0){10}}\put(0,0){\line(1,1){10}}\put(-6,-2){\scriptsize 1}\put(-6,8){\scriptsize 4}\put(13,-2){\scriptsize 2}\put(13,8){\scriptsize 3}
\end{picture} with $G = \langle (1\ 3)(2\ 4)\rangle$, $\{i,j\} = \{1,3\}$ and $\{\a,\a'\} = \{\a_2,\a_4\}$.

\end{document}